
\input amstex
\documentstyle{amsppt}
\hoffset=1cm\voffset=.5cm
\pretolerance=10000
\tolerance=10000
\parindent=20pt
\pageno=1
\parskip=10pt
\baselineskip=20pt
\nologo
\NoBlackBoxes
\magnification=1200
\pageheight {19.2cm}
\pagewidth{13.8cm}
\hcorrection {-10mm}
\magnification=1200
\topmatter
\author F. Guill\'en et V. Navarro Aznar
\endauthor
\title  Un crit\`ere d'extension d'un foncteur  \\
d\'efini sur les sch\'emas lisses
\endtitle
\affil
Departament d\'Algebra i Geometria \\
Facultat de Matem\`atiques
- Universitat de Barcelona\\
 Gran Via, 585 - 08007 Barcelone - Espagne\\
guillen\@cerber.mat.ub.es\\
Avril 95
\endaffil
\thanks Ce travail a \'et\'e partiellement subventionn\'e par le projet DGCYT
PB93-0790 \endthanks
 \endtopmatter
\document
\subheading{Introduction}

Soient $k$ un corps de caract\'eristique z\'ero et $X$ une vari\'et\'e
alg\'ebrique
sur $k$. On sait, d'apr\`es le th\'eor\`eme de r\'esolution des singularit\'es
d'Hironaka (\cite{H1}) , qu'on peut r\'esoudre les singularit\'es de $X$,
c'est-\`a-dire, qu'il existe une vari\'et\'e non singuli\`ere $\widetilde X$ et
un morphisme $f:\widetilde X\longrightarrow X$ qui est birationnel et
propre. En plus, si $Y$ est une sous-vari\'et\'e ferm\'ee de $X$, il existe
une r\'esolution   $f:\widetilde X\longrightarrow X$ telle que $\widetilde
Y=f^{-1}(Y)$ soit un diviseur \`a croisements normaux dans $\widetilde
X$.

On a donn\'e de nombreuses applications de ce th\'eor\`eme
de r\'esolution \`a l'\'etude cohomologique des vari\'et\'es alg\'ebriques et,
en
particulier, il a \'et\'e utilis\'e pour \'etendre certains foncteurs
cohomologiques d\'efinis \`a priori sur une classe de sch\'emas lisses \`a une
classe plus vaste de sch\'emas, des exemples de telles extensions \'etant la
cohomologie de De Rham (\cite{Gr},\cite{Ha}) et la th\'eorie de
Hodge-Deligne (\cite{D1}).

Dans cet article nous prouvons, \`a partir du th\'eor\`eme d'Hironaka et comme
continuation de notre  travail pr\'ec\'edent (\cite{HC}), un crit\`ere
d'extension d'un foncteur d\'efinit sur les sch\'emas s\'epar\'es, de type fini
et lisses sur $k$, crit\`ere qui dans un langage peu pr\'ecis, montre que
si
on a pour ce foncteur la suite exacte habituelle d'un \'eclatement, alors
on peut \'etendre ce foncteur \`a tous les sch\'emas s\'epar\'es et de type
fini
sur $k$. Plus pr\'ecis\'ement, nous prouvons au  \S 2 le r\'esultat
suivant.

 Soient $k$ un corps
de caract\'eristique
z\'ero, $\bold{Sch}(k)$ la cat\'egorie  des sch\'emas
s\'epar\'es et de type fini sur $k$ et $\bold{Sch_{Reg}}$(k) la
sous-cat\'egorie des sch\'emas lisses.
Soient $\bold{A}$ une cat\'egorie ab\'elienne, $\text{C}^{b}(\bold{A})$ la
cat\'egorie de complexes born\'es de $\bold{A}$, et $\text{D}^{b}(\bold{A})$
sa cat\'egorie d\'eriv\'ee.

\proclaim{ Th\'eor\`eme}
Soit
$
G:\bold{Sch_{Reg}}(k)\longrightarrow \text{C}^{b}(\bold{A})
$ un foncteur contravariant
 tel que: (F1) $G(\emptyset)= O$; (F2)
il existe un morphisme naturel $G(X\coprod Y)
\longrightarrow  G(X)
\times G(Y)$ qui est un quasi-isomorphisme;
et (F3) si $i:Y\longrightarrow X$ est une immersion ferm\'ee de
$\bold{Sch_{Reg}}(k)$, $f:\widetilde X \longrightarrow X$ est
l'\'eclatement de $X$ le
long de $Y$ , $\widetilde Y$ est le diviseur exceptionnel, $\widetilde
Y=f^{-1}(Y) $, $j:\widetilde Y\longrightarrow \widetilde X $ d\'enote
l'inclusion et $g:\widetilde Y\longrightarrow Y$  la
restriction de $f$, alors le morphisme $$
G(X) @>i^{*}+f^{*}>> \bold{s}\left(G(Y) \bigoplus_{}^{}G(\widetilde X)
@>g^{*}-j^{*}>> G(\widetilde Y)\right) $$
est un quasi-isomorphisme, o\`u $\bold{s} $ d\'enote le complexe simple ou
total.
Alors,
il existe une unique extension de $G$ \`a un foncteur contravariant
$
G:\bold{Sch}(k)\longrightarrow \text{D}^{b}(\bold{A}), $
qui v\'erifie la propri\'et\'e de descente cohomologique
suivante:
si
$f:\widetilde X \longrightarrow X$ est un morphisme propre de
$k$-sch\'emas, $i:Y\longrightarrow X$ est une immersion ferm\'ee,
$\widetilde Y=f^{-1}(Y) $,  $j:\widetilde Y\longrightarrow \widetilde
X $ d\'enote l'inclusion, $g:\widetilde Y\longrightarrow Y$  la
restriction de $f$, et $f$ induit un isomorphisme $\widetilde X -
\widetilde Y \longrightarrow  X -  Y $,
 alors
le morphisme $$
G(X) @>i^{*}+f^{*}>> \bold{s}\left(G(Y) \bigoplus_{}^{}G(\widetilde X)
@>g^{*}-j^{*}>> G(\widetilde Y)\right) $$
est un quasi-isomorphisme.
\endproclaim

Notons que dans ce r\'esultat nous  consid\'erons les foncteurs
cohomologiques comme \'etant des foncteurs prenant des valeurs
dans une cat\'egorie de complexes, en accord avec les axiomatisations
r\'ecentes des th\'eories cohomologiques (\cite{B} et \cite{Gi}), tandis
que dans
la formulation classique de
 la cohomologie
singuli\`ere et de ses diff\'erentes g\'en\'eralisations: cohomologies de Weil
ou de Bloch-Ogus, par exemple, on ne consid\`ere que le foncteur qu'on
obtient par application du foncteur de cohomologie $H^{*}$. L'actuelle
formulation parait essentielle pour traiter les probl\`emes li\'es \`a la
th\'eorie de la descente cohomologique (\cite{D1}), et, en particulier, le
probl\`eme qui nous occupe.

Or, dans
les applications que nous avons en vue du crit\`ere d'extension:
au complexe filtr\'e de Hodge-De Rham, \`a la th\'eorie des motifs, et \`a
la
th\'eorie d'homotopie rationnelle en th\'eorie de De Rham, la cat\'egorie o\`u
prend ses valeurs le foncteur consider\'e ne provient pas d'une
cat\'egorie
ab\'elienne, c'est par exemple le cas de la cat\'egorie des motifs de Chow
qui est pseudo-ab\'elienne, ou m\^eme d'une cat\'egorie non additive,
comme c'est le cas de la cat\'egorie des alg\`ebres
diff\'erentielles
gradu\'ees
commutatives  (dgc)
 sur $k$, n\'ecessaire pour l'homotopie rationnelle. Ainsi, nous sommes
oblig\'es \`a prouver le crit\`ere d'extension dans une
situation plus g\'en\'erale, non
n\'ecessairement additive, dans laquelle on substitue la
cat\'egorie
des complexes d'une cat\'egorie ab\'elienne par une cat\'egorie v\'erifiant
certaines propri\'et\'es naturelles au contexte et que nous appelons
cat\'egorie de descente.
Le concept de cat\'egorie de descente est une variante des
cat\'egories
triangul\'ees de Verdier et nous  developpons ce formalisme dans le
premier paragraphe. Ce
contexte non additif inclut le cas des alg\`ebres dgc
 sur un corps $k$ de caract\'eristique z\'ero.
C'est le cadre qui nous  permet de donner, au \S 3, une
r\'ealisation en th\'eorie de De Rham alg\'ebrique de
l'homotopie rationnelle des $k$-sch\'emas.

Le crit\`ere d'extension pr\'ec\'edent est aussi v\'erifie dans le cas
des espaces
analytiques, c'est-\`a-dire, si on substitue dans l'enonc\'e la cat\'egorie
des $k$-sch\'emas par la cat\'egorie des espaces analytiques et la
sous-cat\'egorie des sch\'emas lisses par celle des vari\'et\'es complexes, car
on a aussi dans ce contexte des th\'eor\`emes de r\'esolution
(\cite{AH}, aussi \cite{BM}).
Comme application  nous prouvons, au \S 4,
l'existence
du complexe filtr\'e de Hodge-De Rham
pour tout espace
analytique sans recours \`a la th\'eorie de Hodge-Deligne (cf.\cite{DB}
et \cite{HC}(Exp.V)), ce complexe pour une
variet\'e complexe $X$ n'\'etant autre que  le complexe $\Omega ^{*}(X)$
des formes diff\'erentielles
holomorphes muni de la filtration par le degr\'e, $F^{p}\Omega
^{*}(X)=\Omega ^{*\ge p}(X)$.

Finalement, au \S 5,   nous
prouvons l'existence de deux foncteurs $h$ et $h_{c}$ qui
\'etendent la th\'eorie des
motifs de Grothendieck \`a la cat\'egorie des sch\'emas $\bold{Sch}(k)$ et qui
correspondent \`a une  th\'eorie cohomologique sans support et \`a support
compact, respectivement. Bien que ces foncteurs ne soient pas les
foncteurs motiviques qu'on attend, par exemple le foncteur $h$
ne serait que le terme $E_{1}$ de la suite spectrale associ\'ee \`a la
filtration par le poids du "vrai motif mixte",
on en d\'eduit une r\'eponse  affirmative au probl\`eme
pos\'e par Serre (\cite{Se},\S 8)
 sur l'ind\'ependance du motif virtuel d'une
vari\'et\'e alg\'ebrique,
ce probl\`eme avait \'et\'e d\'ej\`a r\'esolu par Gillet
et
Soul\'e au congr\'es Algebraic K-Theory Conference de Paris, juillet 1994,
suivant une autre voie
(\cite {GS}).

Nous sommes reconnaissants \`a C. Soul\'e pour les
discussions que nous avons eues sur ce sujet.

\subheading{1. Cat\'egories de descente}

Dans ce paragraphe nous  introduisons les cat\'egories de
descente, qui sont une variante  des cat\'egories
triangul\'ees de Verdier
(\cite{V}) adapt\'ee \`a la formulation
des th\'eories  cohomologiques non additives.

\noindent
$\bold{(1.1)}$
Soient $\bold{D}$ une cat\'egorie et ${E}$ une classe de morphismes
de $\bold{D}$.
D\'enotons par  $\gamma:\bold{D}\longrightarrow{E}^{-1}\bold{D}$ le
foncteur de localisation, c'est-\`a-dire, le foncteur $\gamma $ transforme
les morphismes de ${E}$ en isomorphismes, et il est universel avec
cette propri\'et\'e (voir \cite{GZ}).
 Nous dirons que ${E}$ est une classe satur\'ee de morphismes
si, pour tout morphisme $f$ de
$\bold{D}$, $f$ est dans ${E}$  si, et seulement si, $\gamma (f)$
est un isomorphisme dans ${E}^{-1}\bold{D}$. Par exemple,
si ${E}$ est la classe de morphismes qui sont des isomorphismes
d'apr\`es l'application d'un foncteur, on voit ais\'ement que ${E}$
est satur\'ee.

\proclaim{(1.2) \, Proposition}
Supposons que ${E}$ est une classe satur\'ee de morphismes de
$\bold{D}$, alors on a
\roster
\item Les isomorphismes de $\bold{D}$ sont dans ${E}$.
\item Si $f:X\longrightarrow Y$ et $g:Y\longrightarrow Z$ sont des
morphismes de $\bold{D}$ et deux des trois morphismes ${f,g,g\circ f}$
sont dans ${E}$, alors le troisi\`eme l'est aussi.
\item
Si $f:X\longrightarrow Y$, $g:Y\longrightarrow Z$ et
$h:Z\longrightarrow T$ sont des morphismes de $\bold{D}$
et $g\circ f$, $h\circ g$ sont dans ${E}$, alors  $g$ est
dans ${E}$.
\endroster

\endproclaim

\demo {Preuve} Les propri\'et\'es (1) et (2) sont \'evidentes.
Pour
prouver (3) nous rappelons que,
dans une cat\'egorie arbitraire,
si $g\circ f$ est un isomorphisme, $f$ a une r\'etraction et $g$ a une
section, et si un morphisme $f$ a une r\'etraction et une section, alors
c'est un isomorphisme.
Ainsi, avec les hypoth\`eses de (3), $g$ est un
isomorphisme dans
${E}^{-1}\bold{D}$, et compte tenu que ${E}$ est satur\'ee, il
en r\'esulte  que $g$
est dans ${E}$.

\enddemo

\noindent
$\bold{(1.3)}$
Nous utiliserons dans ce qui suit
la cat\'egorie d'objets cubiques associ\'ee \`a une cat\'egorie.
On note
$\square_{0}^{+}$ l'ensemble ordonn\'e \`a deux \'el\'ements
$[0,1]$.
Le $n$-cube augment\'e $\square^{+}_{n}$ est le produit cart\'esien
de $n+1$ copies de $\square ^{+}_{0}$. Cet ensemble est muni de
l'ordre lexicographique du produit.
Le sous-ordre de $\square_{n}^{+}$ form\'e par
suppression du premier \'el\'ement $(0,...,0)$ se note $\square_{n}$.
On d\'efinit le $\infty $-cube
augment\'e $\square^{+}_{\infty }$ comme l'ensemble des applications
$\alpha
:\Bbb
N \longrightarrow \square^{+}_{0}$, muni de l'ordre lexicographique.

 Nous
appellerons ordre cubique tout sous-ordre fini $\square$ de
$\square^{+}_{\infty }$ tel que
pour tout couple  $\alpha,\,\beta  $, d'\'el\'ements de
$\square$ on a $\{\gamma \in
\square ^{+}_{\infty };\alpha \le\gamma \le\beta \}\subseteq \square$.
En particulier $\square^{+}_{n}$ et $\square^{}_{n}$ sont des
ordres cubiques.
Les ordres cubiques avec les applications
d'inclusion d\'efinissent une cat\'egorie que nous noterons
$\Pi $ .

Soit maintenant $\bold{A}$ une cat\'egorie. On d\'efinit la cat\'egorie
$(\Pi ,\bold{A})$
dont les objets sont les foncteurs
$X_{\bullet}:\square \longrightarrow \bold{A} $ , d'un ordre
cubique
$\square$, appel\'e le type de $X_{\bullet}$, dans $\bold{A}$, avec
$\square $ variable dans $\Pi $. Si $X_{\bullet}:\square
\longrightarrow \bold{A}$ et
$X'_{\bullet}:\square^{'} \longrightarrow \bold{A}$ sont des objets
de $(\Pi ,\bold{A})$  , tels qu'on a une inclusion $\delta
:\square\longrightarrow \square'$, un morphisme de
$(\Pi ,\bold{A})$
$\tau :X_{\bullet}\longrightarrow
X'_{\bullet}$ est une transformation naturelle de
foncteurs $\tau :X_{\bullet}\longrightarrow X'_{\bullet}\circ \delta $.
 On a un foncteur  $type :
(\Pi ,\bold{A})\longrightarrow \Pi $.
La cat\'egorie $(\Pi ,\bold{A}^{op})$ sera denot\'ee par
$(\Pi ^{op},\bold{A})$.

Si $f:X_{0}\longrightarrow X_{1}$ est un morphisme d'une cat\'egorie
$\bold{A}$, nous noterons $tot(f)$ l'objet de $
(\Pi ,\bold{A})
$ de type
$\square^{+}_{0}$ d\'efini par $X_{\bullet}$  (voir
\cite{HC}(I.1.7)).

 \proclaim{(1.4) D\'efinition} On appelle cat\'egorie
de descente cubique (ou simplement de descente) une
cat\'egorie
$\bold{D}$ munie d'une classe de morphismes ${E}$ et
d'un foncteur covariant
$$
\bold{s}:(\Pi ,\bold{D})\longrightarrow
\bold{D}, $$
que nous appellerons foncteur simple, v\'erifiant les
conditions suivantes:
\roster
\item"(C1)." $\bold{D}$ a un objet final, not\'e $O$.
\item"(C2)." $\bold{D}$ a des produits finis.
\item"(E1)." ${E}$ est une classe satur\'ee de morphismes.
\item"(E2)." Si $f$, $g$ sont dans ${E}$, le produit $f \times
g$ est dans ${E}$.
\item"(S1)."
Pour tout objet $X$ de $\bold{D}$,
il existe  un isomorphisme naturel
dans $\bold{D}$
$$
\bold{s}(tot(X\rightarrow O))\longrightarrow X.
$$
\item"(S2)." Si $X_{\bullet}$ et $Y
_{\bullet}
$ sont des objets cubiques   de
$\bold{D}$ du m\^ eme type,
on a un isomorphisme $\bold{s}(X
_{\bullet}
\times Y_{\bullet}) \cong\bold{s}(X_{\bullet})\times
\bold{s}(Y_{\bullet}).$
\item"(S3)." Pour tout morphisme $f$ de $\bold{D}$,
$f$ est dans $ {E}$ si, et seulement si, le
morphisme $\bold{s}( tot(f))\longrightarrow O$ est dans
${E}$.
\item"(S4)." Si $f_{\bullet}:X
_{\bullet}
\longrightarrow Y_{\bullet} $ est un morphisme
d'objets cubiques de type $\square$, et $f_{\alpha }$ est dans
${E}$, pour
tout $ \alpha \in \square$, alors $\bold{s}(f
_{\bullet}
)$ est dans ${E}$.
\item"(S5)." Si $\square= \square ' \times \square ''$, pour tout
objet cubique $X_{\bullet\bullet}$ de $
\bold{D}$ de type $\square$,  il existe des morphismes  dans ${E}$
$$
\bold{s}_{\square ''}(\bold{s}_{\square'}(X
_{\bullet\bullet}))\longleftarrow
\bold{s}_{\square}(X
_{\bullet
\bullet})\longrightarrow \bold{s}_{\square
'}(\bold{s}_{\square''}(X_{\bullet\bullet}))
$$
qui sont naturels.
\endroster
\endproclaim

 Nous d\'enoterons par $Ho\,\bold{D}$
la cat\'egorie localis\'ee ${E}^{-1}\bold{D}$. Nous appellerons, en
suivant l'usage,  quasi-isomorphismes (dans d'autres
contextes ils
sont appel\'es \'equivalences faibles)  les morphismes de ${E}$, et
acycliques les objets de $\bold{D}$ quasi-isomorphes \`a l'objet final
$O$.

\remark{ Remarque} La donn\'ee d'un foncteur simple
$$
\bold{s}:(\Pi ,\bold{D})\longrightarrow
\bold{D} $$
est \'equivalente
\`a la donn\'ee, pour tout ordre cubique $\square$,
d'un foncteur
$$
\bold{s}_{\square}:\bold{D}^{\square}\longrightarrow
\bold{D} $$
naturel par rapport aux inclusions des ordres cubiques.
\endremark

 \noindent $\bold{(1.4)^{op}}$
Les axiomes envisag\'es pr\'ec\'edemment ne sont pas autoduals
dans la variable $\bold{D}$, et ils sont appropri\'es pour une th\'eorie de
la descente cohomologique, mais
il existe une d\'efinition duale, pour une th\'eorie de la descente
homologique: on dit que $\bold{D}$ est une
cat\'egorie de descente cubique
homologique
si la cat\'egorie oppos\'ee $\bold{D}^{op}$ est de descente cubique.

\noindent
$\bold{(1.5)}$ Si dans la d\'efinition pr\'ec\'edent on substitue les
ordres
cubiques $\square$ par les produits finis de cat\'egories simpliciales
strictes tronqu\'ees inf\'erieurement et sup\'erieurement, on obtient
une notion de cat\'egorie de descente simpliciale. D'apr\`es
\cite{Gu}(2.1.6), il est ais\'e de voir que toute cat\'egorie de descente
simpliciale est aussi une cat\'egorie de descente cubique.

\noindent $\bold{(1.6)\, Exemple.} $
Soient $\bold{A}$ une cat\'egorie ab\'elienne,
$C^{\alpha }(\bold{A})$ la cat\'egorie
des complexes de cocha{\^\i}nes de $\bold{A}$ (o\`u $\alpha =+,b,\emptyset)$,
${E}$ la classe des homologismes
(aussi appel\'es quasi-isomorphismes), et $\bold{s}$
le foncteur simple ordinaire d'un complexe cubique (voir
\cite{HC},(I.6.1)).
Alors il est ais\'e de v\'erifier les conditions de cat\'egorie de
descente pour les donn\'ees pr\'ec\'edentes. Remarquons en
particulier que (S4)
est une cons\'equence de la suite spectrale de la filtration par l'index
cubique, qui est toujours bir\'eguli\`ere, et que
la propri\'et\'e (S5) s'obtient de l'identit\'e
$\bold{s}_{\square}$=$\bold{s}_{\square '} \, \circ \,
\bold{s}_{\square ''}$.

\noindent $\bold{(1.7) \, Exemple.}$
Soient $\bold{A}$ une cat\'egorie additive, $
C^{\alpha }(\bold{A}) $ la
cat\'egorie des complexes de cocha{\^\i}nes de $\bold{A}$ (o\`u $\alpha  =
b, +$, $\emptyset$), ${E}$ la classe des
homotopismes, c'est-\`a-dire, les \'equivalences homotopiques,
et $\bold{s}$ le foncteur simple ordinaire.
Nous prouverons \`a la suite que $
C^{\alpha }(\bold{A})
$ est une cat\'egorie de
descente.

Les propri\'et\'es (C1), (C2), (E2), (S1) et (S2) sont \'evidentes. La
preuve de (S5)
est comme dans le cas d'une cat\'egorie ab\'elienne. La propri\'et\'e (S3)
est un r\'esultat
classique: Un morphisme de complexes est une \'equivalence homotopique
si et seulement si le c\^one du morphisme est contractile.
Finalement, les conditions (E1) et (S4)
r\'esultent de la proposition
suivante.

\proclaim{(1.8) Proposition}
Avec les hypoth\`eses de {\rm (1.7)} on a
\roster
\item La cat\'egorie $\pi\,
C^{\alpha }(\bold{A})
 $, quotient de $C^{\alpha }(\bold{A})$ par la
relation
d'homotopie, est \'equivalente \`a la cat\'egorie localis\'ee
$Ho \,C^{\alpha }(\bold{A})$.
\item La classe ${E}$ est satur\'ee.
\item
Soient $n$ un entier $\ge 0$, $f_{\bullet}:X
_{\bullet}
\longrightarrow Y_{\bullet}$
un morphisme de $\square^{+}_{n}$-complexes, et supposons que $f_{\alpha
}$ soit une \'equivalence homotopique, pour tout $\alpha \in
\square^{+}_{n}$,
alors le morphisme $\bold{s}(f_{\bullet}):\bold{s}(X
_{\bullet}
)\longrightarrow
\bold{s}(Y_{\bullet})$ est  une \'equivalence homotopique.
\endroster
\endproclaim

\demo{Preuve}
La propri\'et\'e (1) est une cons\'equence de \cite{Q}(I.1), lemma 8.
Rappelons l'argument par la commodit\'e du lecteur.  D'abord, on remarque
que si un morphisme $f$ a des inverses homotopiques \`a gauche $g$,
et \`a droite $g'$:
$$
f\circ g' \backsim 1 , \qquad g\circ f \backsim 1,
$$
alors on a $[g]=[g']$,  et donc [f] a un inverse dans $\pi\,
C^{\alpha }(\bold{A})$,
d'o\`u il r\'esulte que les homotopismes sont des isomorphismes dans
$\pi\,C^{\alpha }(\bold{A}) $.

D'autre part,
il existe un foncteur cylindre $Cyl:C^{\alpha }(\bold{A})
\longrightarrow C^{\alpha }(\bold{A})$,
et des homoto\-pismes naturels en $X$,
$$
\delta _{0},\delta _{1}:X\longrightarrow Cyl(X), \quad \sigma
:Cyl(X)\longrightarrow X,
$$
tels que $\sigma \circ \delta _{0}=\sigma \circ \delta _{1}=1_{X}$,
pour tout objet $X$ de $C^{\alpha }(\bold{A})$.

Maintenant, soit
$$
T:C^{\alpha }(\bold{A})\longrightarrow \bold{U}
$$
un foncteur qui transforme  les homotopismes en  isomorphismes. Si
$f,g:X\longrightarrow Y$ sont des morphismes homotopes, on a une
homotopie $H:f\backsim g$, c'est-\`a-dire, un morphisme
$$
H:Cyl(X)\longrightarrow Y
$$
tel que $H\circ \delta _{0}=f$ et $H\circ \delta _{1}=g$, d'o\`u il
r\'esulte
$$T(f)=TH\circ T\delta _{0}=TH \circ (T\sigma )^{-1}
=TH\circ T\delta_{1} =T(g),
$$
il s'ensuit que $T$ factorise uniquement \`a travers de $\pi
\,C^{\alpha }(\bold{A})$, ce qui donne (1).

 La propri\'et\'e (2) r\'esulte imm\'ediatement de (1). En effet,
si un morphisme $f$ de $
C^{\alpha }(\bold{A})
$ est tel que $[f]$ est un isomorphisme
dans la cat\'egorie $\pi       \,
C^{\alpha }(\bold{A})$, et $[g]$ est l'inverse de $[f]$,
alors $g$ est un inverse homotopique de $f$, et donc $f$ est un
homotopisme.

Pour prouver (3) proc\'edons par r\'ecurrence sur $n$. Pour $n=0$ on a
la situation suivante: $X_{\bullet}$ est un morphisme de complexes
$\xi  :X_{0}\longrightarrow X_{1}$ ,  $Y_{\bullet}$ est un morphisme
de complexes $\eta :Y_{0}\longrightarrow Y_{1} $ , et $f
_{\bullet}
$ est un
couple de morphismes $(f_{0},f_{1})$ tel que le diagramme suivant
$$
\CD
X_{0} @>f_{0}>>Y_{0}\\
@V\xi VV @VV \eta V \\
X_{1} @>f_{1}>> Y_{1}
\endCD
$$
est commutatif.
Puisque
la cat\'egorie $K^{b}(\bold{A})$ des complexes \`a homotopie pr\`es est
triangul\'ee (voir \cite{V}(I.2.2)), on a que
 $\bold{s}(f_{\bullet})$
est une \'equivalence homotopique
(voir \cite{V}(I.1.2)).

Supposons maintenant que $n$ soit  $\ge$ 1. On d\'ecompose
$\square^{+}_{n}=\square^{+}_{0} \times \square^{+}_{n-1}$, et on
utilise la propri\'et\'e  (S5).
Il r\'esulte, d'apr\`es l'hypoth\`ese de r\'ecurrence  que
$\bold{s}_{n-1}(f_{1\bullet})$ et $\bold{s}_{n-1}(f_{0\bullet})$ sont
des \'equivalences
homotopiques, et d'apr\`es le cas $n=0$, on en d\'eduit que
$\bold{s}(f_{\bullet})$ est une \'equivalence homotopique.

\enddemo

 \noindent $\bold{ (1.9) \, Exemple.}$  De fa\c con analogue que dans
l'exemple pr\'ec\'edent, on peut v\'erifier
les conditions de cat\'egorie de descente
pour les donn\'ees suivantes: $\bold{Top}_{0}$ la cat\'egorie des
espaces topologiques connexes par arcs et point\'es,
${E}$ la classe des \'equivalences faibles
et $\bold{s}$
le foncteur simple
d\'efini par
$$
\bold{s}X_{\bullet}=\lim_{\to } sc_{P}X_{\bullet},
$$
o\`u $P:\bold{Top}_{0}\longrightarrow \bold{Top}_{0} $ est le foncteur
espace de chemins $P(X,x)=(X,x)^{(I,0)}$ (voir
\cite{HC}(VI.2) pour les notations). En particulier, dans le cas
d'un diagramme tot ($f$) associ\'e
\`a une application continue $f:(X,x)\longrightarrow (Y,y)$, on a
$\bold{s}({\text tot}(f))=(X,x)\times_{(Y,y)}P(Y,y) $, qui est la
fibre homotopique
 ou espace de chemins de l'application $f$.

 \noindent
$\bold{(1.10)}$ Le lemme
suivant  sera  utilis\'e pour v\'erifier les conditions
de cat\'egorie de descente dans les exemples que nous donnerons dans les
applications. La preuve du lemme  est imm\'ediate.

\proclaim{ Lemme} Soit $\bold{D}$ une cat\'egorie munie d'un foncteur
$$
\bold{s}:(\Pi ^{op},\bold{D})\longrightarrow \bold{D}
$$
qui v\'erifie {\rm(C1), (C2), (S1)} et {\rm (S2)}.

Supposons qu'il existe une cat\'egorie de descente $\bold{D'}$
et un foncteur
$$
\phi :\bold{D}\longrightarrow \bold{D'}
$$
tel que
\roster
\item $\phi (X\times Y) = \phi(X)\times \phi (Y) $, et $\phi
(O)= O'$.
\item Pour tout objet cubique $X.$ de $\bold{D}$ il existe un morphisme
naturel $$\phi (\bold{s}(X.))\longrightarrow \bold{s'}(\phi(X.))$$ dans
${E}$.
\endroster
Si on d\'efinit
$$
E=\{f\in Mor\bold{D};\phi (f)\in E'\}
$$
alors $\bold{D}$ est une cat\'egorie de descente.
\endproclaim

\noindent$\bold{(1.11)}$
Soit $\bold{D}$ est une cat\'egorie de descente,
munie
d'une classe distingu\'ee de morphismes  et d'un foncteur
simple, not\'es ${E}$ et
$\bold{s}$ respectivement.
Pour toute petite cat\'egorie
$\bold{I}$, la cat\'egorie
$\bold{D}^{\bold{I}}$
des foncteurs $\bold{I}\rightarrow \bold{D}$
 est aussi une
cat\'egorie de descente, si on d\'efinit  ${E}$ et $\bold{s}$ par
composantes. De la m\^ eme fa\c con, $(\Pi ,\bold{D})$ est une
cat\'egorie de descente. Il en r\'esulte, d'apr\`es
la propri\'et\'e (S4), que pour tout ordre cubique $\square$ le foncteur
$\bold{s}_{\square}$ induit un foncteur
$$
\bold{s}_{\square}:Ho\,(\bold{D}^{\square})\longrightarrow
Ho\,\bold{D}
$$
naturel par rapport aux inclusions des ordres cubiques et donc un
foncteur $$
\bold{s}:Ho\,(\Pi ,\bold{D})\longrightarrow
Ho \,\bold{D}.
$$

\noindent
$\bold{(1.12)}$
Rappelons (\cite{HC}(I.1.9)) qu'une cat\'egorie ordonnable finie est
une cat\'egorie avec un nombre fini de morphismes, sans endomorphismes
diff\'erents des identit\'es et telle que si on munit l'ensemble
des objets de la relation de pr\'eordre: "$i \le j$ si et seulement si
$Hom\,(i,j)$ est non vide" , alors ce pr\'eordre est un ordre.
Une sous-cat\'egorie d'une cat\'egorie ordonnable finie est,
\'evidemment,
une cat\'egorie de m\^ eme esp\`ece, et  ces cat\'egories, avec les
foncteurs d'inclusion, forment une cat\'egorie que nous noterons $\Phi $.
En outre, on
remarque que tout ordre est une cat\'egorie ordonnable, en particulier
tout ordre cubique $\square$
est une cat\'egorie ordonnable finie, et on a un foncteur d'inclusion
$\Pi \longrightarrow \Phi $.

Pour toute cat\'egorie $\bold{A}$ nous noterons, comme dans (1.3),
$(\Phi ,\bold{A})$
la cat\'egorie form\'ee par les diagrammes d'objets de
$\bold{A}$ de type $\bold{I}$,
avec $\bold{I}$ une cat\'egorie ordonnable finie variable. On a un
foncteur $type:(\Phi ,\bold{A})\longrightarrow \Phi $.
Nous noterons $(\Phi ^{op},\bold{A})$ la cat\'egorie
$(\Phi ,\bold{A}^{op})$.

Soit $\bold{D}$ une cat\'egorie de descente.
D'une fa\c con analogue \`a (1.11) on a une
structure de cat\'egorie de descente sur
$(\Phi ,\bold{D})$ telle que le foncteur $type
:(\Phi ,\bold{D})\longrightarrow \Phi $  factorise par la
                                        cat\'egorie
localis\'ee et on a donc un foncteur  $$type :Ho \,
(\Phi ,\bold{D})\longrightarrow \Phi .$$

\subheading{2. Un crit\`ere d'extension d'un
foncteur}

 Dans ce paragraphe nous donnons un crit\`ere d'extension
d'un foncteur
d\'efini sur les sch\'emas  s\'epar\'es, de type fini et lisses sur un
corps $k$. Ce crit\`ere \'etant bas\'e  essentiellement sur la
r\'esolution
de singularit\'es, comme elle est developp\'ee dans la th\'eorie
des hyperr\'esolutions cubiques de \cite{HC}, nous nous
placerons dans un contexte o\`u cette
th\'eorie soit disponible, par exemple et \`a l'heure actuelle, en supposant
le corps $k$ de caract\'eristique z\'ero.

 \noindent $\bold{(2.1)}$
Soient $k$ un corps de caract\'eristique z\'ero et  $\bold{Sch}(k)$
la cat\'egorie  des sch\'emas
s\'epar\'es et de type fini sur $k$, que nous appellerons
simplement sch\'emas.
Nous dirons qu'un diagramme  de $\bold{Sch}(k)$
$$
\CD
\widetilde Y @> j >> \widetilde X \\
@V g VV @VV f V\\
Y @> i >> X
\endCD
$$
est un isomorphisme relatif propre, si le diagramme est cart\'esien,
$i$ est une  immersion
ferm\'ee, et $f$ est propre
 et  induit un isomorphisme $\widetilde X -\widetilde Y
\longrightarrow  X -  Y $.

Notons $\bold{W}$
la cat\'egorie $\bold{Sch}(k)$ et $\bold{M}$ la sous-cat\'egorie pleine
de $\bold{W}$ dont les objets sont des sch\'emas lisses.
Soit $\bold{D}$ une cat\'egorie de descente, munie
d'une classe distingu\'ee de morphismes  et d'un foncteur
simple, not\'es ${E}$ et
$\bold{s}$ respectivement.
Si on a un foncteur contravariant
$$
G:\bold{M}\longrightarrow \bold{D},
$$
alors, pour toute cat\'egorie $\bold{I}$, $G$ induit un foncteur
contravariant
 $$
G:\bold{M}^{\bold{I}^{op}}\longrightarrow
\left( \bold{D}^{\bold{I}} \right),$$
par $$
G(X
_{\bullet}
)_{\alpha }=G(X_{\alpha }),\quad  \text{si} \quad \alpha \in
\bold{I} $$
qui est  naturel en $\bold{I}$,  et donc $G$
d\'efinit un foncteur contravariant
$$ G:(\Phi ^{op},\bold{M})\longrightarrow
Ho \left(\Phi ,\bold{D}\right),$$
qui commute avec les foncteurs $type $ et $tot$ (voir \cite{HC}(I.6)).

Le probl\`eme qui se pose alors est de donner des conditions naturelles
pour qu'on puisse \'etendre ce foncteur \`a toute la cat\'egorie $\bold{W}$,
et dans ce sens nous donnons le crit\`ere suivant.

\proclaim{(2.2) Th\'eor\`eme} Soit
$$
G:(\Phi ^{op},\bold{M})\longrightarrow  Ho\,\left(
\Phi ,\bold{D}\right)
$$
un foncteur contravariant qui commute avec les foncteurs
{\rm type} et {\rm tot}, et
 tel que: \roster
\item "(F1)"
$G(\emptyset)= O$.
\item"(F2)" Il existe un morphisme naturel $G(X\coprod Y)
\longrightarrow  G(X)
\times G(Y)$ qui est un quasi-isomorphisme,
\item "(F3)" Si $i:Y\longrightarrow X$ est une immersion ferm\'ee de
$\bold{M}$, $f:\widetilde X \longrightarrow X$ l'\'eclatement de $X$ le
long de $Y$ et $X_{\bullet}$ est le diagramme cart\'esien
 $$
\CD
\widetilde Y @> j >> \widetilde X \\
@V g VV @VV f V\\
Y @> i >> X
\endCD \
$$
construit \`a partir de $f$ et $i$,
l'objet $\bold{s}G(X_{\bullet})$ de $Ho \,\bold{D}$ est
acyclique.
\endroster
Alors,
il existe une extension de $G$ \`a un foncteur contravariant,
$$
G': (\Phi ^{op},\bold{W})\longrightarrow Ho
\left( \Phi ,\bold{D}\right), $$
qui commute avec les foncteurs $type$ et $tot$, et
v\'erifie la propri\'et\'e de descente suivante:
\roster
\item "(D)" Si le diagramme $X_{\bullet}$ d\'efini par
$$
\CD
\widetilde Y @> j >> \widetilde X \\
@V g VV @VV f V\\
Y @> i >> X
\endCD
$$
est un isomorphisme relatif propre
de $\bold{W}$, alors l'objet
$\bold{s}G'(X_{\bullet})$
de $Ho \,\bold{D}$
est acyclique.
\endroster
\endproclaim

\demo{(2.3) Preuve de (2.2) }

 Soit $\bold{I}$ une cat\'egorie ordonnable finie.
 Si $X:\bold{I}^{op}\longrightarrow \bold{W}$ est un objet de $(\Phi
,\bold{W})$ de type $\bold{I}$ on d\'efinit l'objet $G'(X)$ de $Ho \,
(\bold{D^{I}})$ par $G'(X):=\bold{s}G(X_{\bullet})$, o\`u
$X_{\bullet}\longrightarrow X$ est une hyperr\'esolution cubique de $X$.
Alors, si
$\widetilde X_{\bullet *}\longrightarrow X_{\bullet}$ est une
hyperr\'esolution cubique d'un
$\bold{I}$-sch\'ema $X_{\bullet}$ ,  pour tout  $\alpha \in
\bold{I}$, $\widetilde X_{\alpha *}\longrightarrow X_{\alpha }$ est une
hyperr\'esolution cubique de $X_{\alpha }$ (voir \cite{HC}(I.2.14)), donc
$G'$ est naturel par rapport \`a
$\bold{I}$ et il commute avec les foncteurs $tot$ et $type$.
En outre,
le th\'eor\`eme pour $\bold{I}$
est une cons\'equence du cas o\`u
$\bold{I}$ est la cat\'egorie ponctuelle, cas auquel on se ram\`ene.

Dans ce qui suit nous noterons
$X^{+}_{\bullet}$ le
diagramme total associ\'e \`a
une augmentation
 $X_{\bullet}\longrightarrow X$
d'un diagramme cubique $X_{\bullet}$.

Pour v\'erifier que $G'$ est bien d\'efini, d'apr\`es \cite{HC}(I.3.10) il
suffit de prouver que
si $X'_{\bullet}\longrightarrow X$ et
$X_{\bullet}''\longrightarrow X$ sont des
hyperr\'esolutions cubiques de $X$, et
$
X'_{\bullet}\longrightarrow  X''_{\bullet}
$
est un morphisme de diagrammes de sch\'emas  sur $X$,
alors
$\bold{s} G(X''_{\bullet})\longrightarrow \bold{s} G(X'_{\bullet})$
est un quasi-isomorphisme.

Nous ferons la d\'emonstration du th\'eor\`eme (2.2)
par r\'ecurrence sur la dimension  de $X$ en
prouvant succesivement les assertions $(A_{n}),\, (B_{n})$ et $(C_{n}),
\, n\ge0, $
 suivantes:

\enddemo
\proclaim{$(A_{n})$}  Si
$X_{\bullet}\longrightarrow X$
est une hyperr\'esolution cubique d'un $k$-sch\'ema lisse $X$ de
dimension
$\le n$, alors $\bold{s}G(X^{+}_{\bullet})$ est acyclique. \endproclaim

\proclaim{$(B_{n})$}
Si le diagramme $X^{+}_{\bullet}$ d\'efini par
$$
\CD
\widetilde Y @>>> \widetilde X \\
@VVV @VVV \\
Y @>>> X,
\endCD
\tag "(D1)"
$$
est une 2-r\'esolution d'un $k$-sch\'ema $X$ de dimension  $ \le n$, et
 $Z^{+}_{\bullet * }\longrightarrow X^{+}_{\bullet }$
est une
hyperr\'esolution cubique 1-it\'er\'ee de $X^{+}_{\bullet}$, alors
$\bold{s}G(Z^{+}_{\bullet *})$ est acyclique.
 \endproclaim

\proclaim{ $(C_{n})$}
 Soit $X$ un $k$-sch\'ema de dimension $\le n$.
Si $X'_{\bullet}\longrightarrow X$ et
$X_{\bullet}''\longrightarrow X$ sont des
hyperr\'esolutions cubiques de $X$, et
$X'_{\bullet}\longrightarrow  X''_{\bullet}$
est un morphisme de diagrammes de sch\'emas sur $X$, alors
$\bold{s} G(X''_{\bullet})\longrightarrow \bold{s} G(X'_{\bullet})$
est un quasi-isomorphisme, en particulier $G'$ est bien d\'efini sur les
sch\'emas de dimension $\le n$ et il v\'erifie la condition (D).
\endproclaim

Nous prouverons $(A_{n}),\, (B_{n})$ et $ (C_{n})$ suivant le sch\'ema:
$(C_{n-1})\Longrightarrow(A_{n}), \, (C_{n-1}) +(A_{n})\Longrightarrow
(B_{n}),\,  (C_{n-1})+(B_{n})\Longrightarrow (C_{n}) $. L'assertion
$(A_{0})$ \'etant triviale, le th\'eor\`eme r\'esulte de $(C_{n})$, pour tout
$n$.

\demo{(2.3.1) Preuve de l'implication $(C_{n-1})\Longrightarrow
(A_{n})$}

Soit $X_{\bullet  }\longrightarrow X$
une hyperr\'esolution cubique $m$-it\'er\'ee.
La preuve de $(A_{n})$ se r\'eduit au cas $m=1$ de la fa\c con suivante.  Il
existe une hyperr\'esolution
cubique $m-1$ it\'er\'ee $X'_{\bullet}\longrightarrow X$, telle que
$X_{\bullet }=(Z_{\bullet *}\longrightarrow X'_{\bullet})$ est une
hyperr\'esolution cubique $1$-it\'er\'ee de $X'_{\bullet}$. Alors
pour tout $\alpha $, $Z_{\alpha  *}\longrightarrow X'_{\alpha }$ est une
hyperr\'esolution cubique $1$-it\'er\'ee de $X'_{\alpha }$. D'apr\`es le cas
$m=1$ il r\'esulte que $\bold{s}(G(Z^{+}_{\alpha *}))$ est acyclique pour
tout $\alpha $ , et, compte tenu de (S5) et (S3), ceci entra{\^\i}ne que
$\bold{s}(
G( X^{+}_{\bullet }) )$ est acyclique.

Maintenant supposons que $X_{\bullet}\longrightarrow X$ est une
hyperr\'esolution
cubique 1-it\'er\'ee de $X$. Si $X_{\bullet}$ est un $\square_{r
}$-sch\'ema, nous proc\'edons par r\'ecurrence sur $r $ .
L'assertion $(A_{n})$  est triviale si $r =0$, donc nous supposerons
que
$r > 0$ .

Par (F1), (F2), (E2) et (S2) on se ram\`ene ais\'ement au cas $X$
irr\'eductible car, $X$ \'etant lisse, les composantes irr\'eductibles sont
les composantes connexes.
Le sch\'ema cubique $X_{\bullet}\longrightarrow X$ \'etant une
hyperr\'esolution
cubique 1-it\'er\'ee de $X$, par (S3), (S4) et (S5) et l'hypoth\`ese de
r\'ecurrence sur $r$ et $(C_{n-1}) $ on se ram\`ene
au cas o\`u  $X^{+}_{\bullet}$ est une  2-r\'esolution  de $X$
$$
\CD
\widetilde Y @>>> \widetilde X \\
@VVV @VVV \\
Y @>>> X
\endCD
$$
Si $Y=X$, $(A_{n})$ r\'esulte de l'hypoth\`ese de r\'ecurrence sur
$r $. Supposons
donc que $Y$ est un sous-sch\'ema propre de $X$.
Consid\'erons  d'abord le cas particulier o\`u le morphisme $\widetilde
X\longrightarrow
X$ est la  composition d'une suite d'\'eclatements avec  des
centres
lisses. Nous prouverons $(A_{n})$  par r\'ecurrence sur le nombre $e$
d'\'eclatements.

Lorsque $e=1$, on a un diagramme commutatif
$$
\CD
\widetilde Z @>>>\widetilde Y @>>> \widetilde X \\
@VVV @VVV @VVV\\
Z@>>> Y @>>> X
\endCD
$$
o\`u $Z$ est un $k$-sch\'ema lisse contenu dans $Y$ qui est le centre de
l'\'eclatement $\widetilde X\longrightarrow
X$, et par cons\'equent on a des isomorphismes $\widetilde X -
\widetilde Z \longrightarrow X - Z $, et
$\widetilde Y - \widetilde Z
\longrightarrow
Y - Z $.
Dans ce qui suit nous noterons
$G(X,Y)=\bold{s}(G(X)\longrightarrow G(Y))$  si on a un morphisme
$Y\longrightarrow X$.
$(C_{n-1})$  entra{\^\i}ne que le morphisme
$
G'(Y,Z)\longrightarrow
G'(\widetilde Y,\widetilde Z)
$
est un quasi-isomorphisme, et par (F3), il r\'esulte que
$
G'(X,Z)\longrightarrow G'(\widetilde X,\widetilde Z)
$
est aussi un quasi-isomorphisme. D'apr\`es (1.2.2) et (S3),
 on en conclut
 que $
\bold{s}G(X^{+}_{\bullet})
$
est acyclique,
ce qui prouve  $(A_{n})$ dans ce cas.

Supposons maintenant $e>1$. Alors on a un
diagramme commutatif
$$
\CD
\widetilde X @>>> X' @>>>  X\\
@AAA @AAA @AAA\\
\widetilde Y @>>>  Y' @>>>Y
\endCD
$$
o\`u $\widetilde X\longrightarrow X'$ s'obtient par \'eclatement de $X'$
avec un centre lisse contenu dans $Y$, et $X'\longrightarrow X$
s'obtient par
composition d'une suite de $e-1$ \'eclatements avec des centres lisses.

D'apr\`es l'hypoth\`ese de r\'ecurrence sur $e$, les morphismes
$$
G'(X,Y)\longrightarrow
G'(X',Y'), \text{ et } G(X',Y')\longrightarrow
G'(\widetilde X,\widetilde Y)
$$
sont des quasi-isomorphismes, et, compte tenu de (1.2.2) et (S3), il
en r\'esulte que $\bold{s}G(X^{+}_{\bullet})$ est acyclique, ce qui prouve
$(A_{n})$ dans ce cas.

Dans le cas g\'en\'eral,
d'apr\`es le lemme de Chow
d'Hironaka (\cite {H1}), \'etant
donn\'e le morphisme $\widetilde X\longrightarrow X$,
il existe
un
morphisme $X'\longrightarrow X$
qui est
la   composition d'une suite
d'\'eclatements avec des centres lisses,
et qui domine  $
 \widetilde X
$. Il existe aussi un
morphisme $\widetilde X'\longrightarrow \widetilde X$ qui est
la composition d'une suite d'\'eclatements avec des centres lisses et qui
domine  $X'$ .
 On a donc un diagramme
commutatif
$$
\CD
\widetilde X' @>>> X' @>>> \widetilde X @>>> X \\
@AAA @AAA @AAA @AAA\\
\widetilde Y' @>>> Y' @>>> \widetilde Y @>>> Y \\
\endCD
$$
tel que les morphismes
$$
\widetilde X '- \widetilde Y
' \longrightarrow X' - Y'\longrightarrow
\widetilde X - \widetilde Y
 \longrightarrow X - Y
$$
sont des isomorphismes. En appliquant $(A_{n})$ dans le cas particulier
d\'emontr\'e auparavant, il r\'esulte que les morphismes $$
G'(X,Y)\longrightarrow
G(X',Y')
\quad , \quad
G(\widetilde X,\widetilde
Y)\longrightarrow
G'(\widetilde X',
\widetilde Y')
$$
sont des quasi-isomorphismes,
et, d'apr\`es (1.2.3) et (S3) appliqu\'es \`a la composition
$$
G'(X,Y)\longrightarrow G'(\widetilde X,\widetilde Y)\longrightarrow
G'(X',Y')\longrightarrow G'(\widetilde X',\widetilde Y'), $$
on
conclut que
$
G'(X,Y)\longrightarrow
G'(\widetilde X,
\widetilde Y)
$
est aussi un quasi-isomorphisme,
ce qui prouve $(A_{n})$.
\enddemo

\demo{(2.3.2) Preuve de l'implication $(C_{n-1}) +
(A_{n})\Longrightarrow B_{n}$}

Par r\'ecurrence noeth\'erienne, nous pouvons supposer que $X$ est
irr\'eductible et que $Y$ est un sous-sch\'ema ferm\'e propre
de $X$. Alors, en proc\'edant par r\'ecurrence sur l'ordre de
l'hyperr\'esolution cubique $Z_{\bullet *}^{+}$
de $X_{\bullet}^{+}$, on se ram\`ene
comme pr\'ec\'edemment  au cas o\`u $Z_{\bullet *}^{++}$ est une
2-r\'esolution.

Comme $Z^{++}_{\bullet *}$ est une 2-r\'esolution de (D1),
en utilisant l'index * pour la 2-r\'esolution, $Z_{\bullet *}^{++}$
est un diagramme
$$
\CD
\widetilde Y_{*}^{+} @>>> \widetilde X^{+}_{*}\\
@VVV @VVV \\
Y^{+}_{ *} @>>> X^{+}_{*},
\endCD
$$
o\`u, d'apr\`es \cite{HC},(I.2.8),
$ Y_{ *}^{+},   \widetilde Y^{+}_{*},  \widetilde X^{+}_{*}$
et $X^{+}_{*}$ , sont des  2-r\'esolutions de $Y, \widetilde Y,
 \widetilde X$ et $X$ respectivement, qui sont les files du diagramme
d\'epli\'e
$$
\CD
Y  @<<< Y_{01} @ <<< Y_{11} @>>>Y_{10}  @>>> Y \\
@AAA @AAA @AAA @AAA @AAA   \\
\widetilde Y  @<<< \widetilde Y_{01}  @ <<< \widetilde Y_{11} @>>>
                      \widetilde Y_{10}  @>>> \widetilde Y \\
@VVV @VVV @VVV @VVV @VVV     \\
\widetilde X  @<<<\widetilde X_{01} @ <<<  \widetilde X_{11} @>>>
                           \widetilde X_{10} @>>> \widetilde X  \\
@VVV @VVV @VVV @VVV @VVV       \\
 X @<<< X_{01} @ <<< X_{11}  @>>>X_{10}  @>>> X \\
@AAA @AAA @AAA @AAA @AAA \\
Y  @<<< Y_{01} @ <<< Y_{11} @>>>Y_{10}  @>>> Y
\endCD
$$
o\`u nous identifions la premi\`ere et la derni\`ere colonnes et proc\'edons de
m\^eme pour les files. Nous remarquons que les sch\'emas avec l'index
$_{10}$ sont lisses.

On peut supposer qu'on a  $\widetilde Y \subset \widetilde X_{01}$
et $Y \subset X_{01}$.
En effet, dans le cas contraire nous pouvons d\'efinir $X'_{01}$
et $\widetilde X_{01}'$ par
$$
X'_{01}:= X_{01} \bigcup Y \, , \, \widetilde X_{01}':= \widetilde
X_{01} \bigcup \widetilde Y,
 $$
et $X'_{11}$, $\widetilde X_{11}'$ de telle fa\c con que les diagrammes
$$
\CD
X_{11}' @>>>  X_{10}  \\
@VVV @VVV \\
X_{01}' @>>> X
\endCD
\foldedtext{\qquad et \qquad $\foldedtext{
$$
\CD
\widetilde X_{11}' @>>>  \widetilde X_{11} \\
@VVV @VVV \\
\widetilde X_{01}' @>>> X
\endCD
$$}$}
$$
soient cart\'esiens. Alors cette substitution ne modifie pas la situation,
puisque la diff\'erence est concentr\'ee dans les diagrammes
$$
\CD
X_{11}' @>>> \widetilde X_{11} \\
@VVV @VVV \\
X_{01} @>>> X_{01}'
\endCD
\foldedtext{\qquad et \qquad $\foldedtext{
$$
\CD
\widetilde X_{11} @>>>  \widetilde X_{11}' \\
@VVV @VVV \\
\widetilde X_{01} @>>> \widetilde X_{01}'
\endCD $$}$}
$$
et, par l'hypoth\`ese de r\'ecurrence ($C_{n-1}$), ces diagrammes induisent
des objets
acycliques, d'apr\`es l'application de $\bold{s}G'$, puisqu'ils sont
cart\'esiens et des isomorphismes relatifs propres.

Supposons  donc  $\widetilde Y \subset
\widetilde X_{01}$ et $Y \subset X_{01}$.
Alors, si  on d\'efinit $\widetilde X'_{01}$ et
$\widetilde X'_{11}$ comme des produits fibr\'es
par les diagrammes cart\'esiens
$$
\CD
\widetilde X_{01}' @>>> \widetilde X \\
@VVV @VVV \\
 X_{01} @>>>  X \\
\endCD
\foldedtext{ \qquad et  \qquad $\foldedtext{
$$
\CD
\widetilde X_{11}' @>>> \widetilde X \\
@VVV @VVV \\
 X_{11} @>>>  X, \\
\endCD
$$}$} $$
on obtient les morphismes de diagrammes cubiques
$$
\CD
\widetilde  Y @>>> \widetilde X _{01}@<<< \widetilde X_{11} @>>>
\widetilde X_{10}@>>>\widetilde Y\\ @VVV@VVV @VVV @VVV @VVV \\
\widetilde  Y @>>> \widetilde X _{01}'@<<< \widetilde X_{11}' @>>>
\widetilde X_{10}@>>>\widetilde Y\\ @VVV@VVV @VVV @VVV @VVV \\
Y @>>>  X _{01}@<<<  X_{11} @>>>  X_{10}@>>> Y \, ,\\
\endCD
\tag "(D2)"
$$
o\`u nous identifions la premi\`ere et la derni\`ere colonnes.

Le foncteur $G'$ appliqu\'e au diagramme cubique d\'efini par le rectangle
inf\'erieur donne un objet cubique acyclique. En effet, le diagramme
$$
\CD
\widetilde Y @>>> \widetilde X _{01}'\\
@VVV @VVV \\
Y @>>> X_{01}
\endCD
$$
est un isomorphisme relatif d'apr\`es la d\'efinition de
$\widetilde X_{01}'$.
Le diagramme $$
\CD
\widetilde X_{11}' @>>> \widetilde X _{10}\\
@VVV @VVV \\
 X_{11} @>>>  X _{10}\\
\endCD
$$
est aussi un isomorphisme r\'elatif, car $Y \subset X_{01}$ et le
diagramme est cart\'esien.

Finalement, l'objet cubique d\'efini par le rectangle sup\'erieur de (D2)
induit un objet acyclique. En effet,
le diagramme
$$
\CD
\widetilde X_{11} @>>>\widetilde X' _{11}\\
@VVV @VVV \\
\widetilde  X_{01} @>>>  \widetilde X _{01}'\\
\endCD
$$
est un isomorphisme relatif, car il est cart\'esien et
$\widetilde Y\subset\widetilde X_{01}  $ ,
 donc le foncteur $G'$
appliqu\'e a ce diagramme donne un objet acyclique.

Il r\'esulte que le diagramme cubique d\'efini par le rectangle ext\'erieur de
(D2),
$$
\CD
\widetilde  Y @>>> \widetilde X _{01}@<<< \widetilde X_{11} @>>>
\widetilde X_{10}@>>>\widetilde Y\\ @VVV@VVV @VVV @VVV @VVV \\
Y @>>>  X _{01}@<<<  X_{11} @>>>  X_{10}@>>> Y \, ,\\
\endCD
$$
induit un objet acyclique. Compte tenu que $\widetilde Y^{+}_{*}$ et
$Y^{+}_{*}$ sont des 2-r\'esolutions de $\widetilde Y$ et $Y$
respectivement, d'apr\`es $(C_{n-1})$ on obtient l'acyclicit\'e de l'objet
induit par le diagramme
 $$
\CD
Y_{01} @ <<< Y_{11} @>>>Y_{10} \\
@AAA @AAA @AAA   \\
\widetilde Y_{01}  @ <<< \widetilde Y_{11} @>>> \widetilde Y_{10} \\
@VVV @VVV @VVV     \\
\widetilde X_{01}  @ <<< \widetilde X_{11} @>>> \widetilde X_{10} \\
@VVV @VVV @VVV       \\
X_{01} @ <<< X_{11} @>>>X_{10}  \\
 @AAA @AAA @AAA \\
Y_{01} @ <<< Y_{11} @>>>Y_{10} \,.
\endCD
$$
C'est objet est  $G(Z^{+}_{\bullet *})$,
donc il en r\'esulte $(B_{n})$.
\enddemo

\demo{(2.3.3)  Preuve de l'implication
$(C_{n-1})+(A_{n})+(B_{n})\Longrightarrow (C_{n}) $}

Montrons que $G'$ est bien d\'efini. Soient
 $X'_{\bullet}\longrightarrow X$ et
$X_{\bullet}''\longrightarrow X$ des
hyperr\'esolutions cubiques de $X$, telles qu'on a un diagramme
commutatif $$
\CD
X'_{\bullet}@.  @>>> @. X''_{\bullet}\\
@. \searrow @.@.  @VVV \\
@.@.@. X
\endCD
$$
alors nous allons v\'erifier que
$\bold{s} G(X''_{\bullet})\longrightarrow \bold{s} G(X'_{\bullet})$
est un quasi-isomorphisme.
Consid\'erons une hyperr\'esolution
cubique du diagramme pr\'ec\'edent,
$$
\CD
Z''_{\bullet *}@<<< Z'_{\bullet *} @>>> Z_{*} @<<< Z''_{\bullet *}\\
@VVV @VVV @VVV @VVV\\
X_{\bullet }'' @<<< X_{\bullet}' @>>> X @<<< X_{\bullet}''
\endCD
$$
Puisque $Z_{\bullet *}'\longrightarrow X_{\bullet}'$
est une hyperr\'esolution cubique d'un sch\'ema cubique lisse, il en
r\'esulte de $(A_{n})$ que $
\bold{s} G(X_{\bullet}')\longrightarrow \bold{s} G(Z_{\bullet *}')
$
est un quasi-isomorphisme et, de fa\c con  analogue,
$\bold{s} G(X_{\bullet}'')\longrightarrow \bold{s} G(Z_{\bullet *}'')$
est un quasi-isomorphisme.
Prouvons que
$\bold{s}G(Z_{\bullet *}')\longrightarrow \bold {s} G(Z_{ *})$
et
$\bold{s}G(Z_{\bullet *}'')\longrightarrow \bold {s} G(Z_{ *})$
sont
aussi des quasi-isomorphismes.  En effet,
supposons que $X^{'+}_{\bullet}$ est un diagramme
$$
\CD
\widetilde Y_{\bullet} @>>> \widetilde X \\
@VVV @VVV \\
Y_{\bullet} @>>> X
\endCD
$$
o\`u $\widetilde Y_{\bullet}\longrightarrow \widetilde Y$ et
$Y_{\bullet}\longrightarrow Y$ sont des hyperr\'esolutions de sous-sch\'emas
ferm\'es $\widetilde Y$ et $Y$ de $\widetilde X$ et $X$ respectivement, et
$$
\CD
\widetilde Y @>>> \widetilde X \\
@VVV @VVV \\
Y @>>> X
\endCD
$$
est une 2-r\'esolution de $X$.
Soit $Z^{++}_{\bullet *}$ le diagramme total  $Z^{'+}_{\bullet
*}\longrightarrow Z^{+}_{*}$, qu'on \'ecrit
$$
\CD
\widetilde Y_{\bullet *}^{+} @>>> \widetilde X^{+}_{* }\\
@VVV @VVV \\
Y^{+}_{\bullet *} @>>> X^{+}_{*},
\endCD
$$
o\`u, d'apr\`es \cite{HC},(I.2.14),
$\widetilde Y^{+}_{*}, Y_{ *}^{+}, \widetilde X^{+}_{*}$
et $X^{+}_{*}$ , sont des  hyperr\'esolutions de $\widetilde Y,
Y, \widetilde X$ et $X$ respectivement. On d\'eduit de l'hypoth\`ese de
r\'ecurrence que
$$G'(\widetilde Y)\longrightarrow \bold{s}G(\widetilde
Y_{\bullet })\longrightarrow \bold{s}G(\widetilde Y_{\bullet *})$$  sont
des quasi-isomorphismes et, de fa\c con analogue,
$$G'(Y)\longrightarrow \bold{s} G(Y_{\bullet })\longrightarrow
\bold{s}G( Y_{\bullet
*})$$
sont des quasi-isomorphismes.  On  d\'eduit de $(B_{n})$ que
$\bold{s}G(Z^{+}_{\bullet *})$ est acyclique, et donc
$\bold{s}G(Z_{\bullet*})\longrightarrow \bold{s}G(Z'_{\bullet *})$ est
un quasi-isomorphisme.
De fa\c con analogue, le morphisme
$
\bold{s}G( Z_{\bullet *}'')\longrightarrow \bold {s}
G(Z_{*}) $
est un quasi-isomorphisme,
d'o\`u on d\'eduit finalement que
$$
\bold{s} G(X_{\bullet}'')\longrightarrow \bold{s} G(X_{\bullet}')
$$
est un quasi-isomorphisme.

Il reste \`a v\'erifier la propri\'et\'e de descente (D).
Par r\'ecurrence noeth\'erienne
on peut supposer que $X$ est
irr\'eductible. Le cas $Y=X$ r\'esulte de l'ind\'ependance de $G'$. Ainsi
on
se ram\`ene au cas o\`u $dim \, Y, dim \, Y' < n$. Dans ce cas, la preuve
est analogue \`a celle de l'assertion $(B_{n})$, en utilisant l'hypoth\`ese
de r\'ecurrence $(C_{n-1})$.

\enddemo

\proclaim{(2.4) Corollaire} Soient
$$
F,G: (\Phi ^{op},\bold{M})\longrightarrow  Ho
\, (\Phi ,\bold{D}) $$
des foncteurs v\'erifiant les conditions {\rm (F1), (F2)} et {\rm (F3)}
du th\'eor\`eme {\rm (2.2)}, et
$$
\tau  :F\longrightarrow G
$$
une transformation naturelle de foncteurs. Si $F'$ et $G'$
sont des extensions de $F$ et $G$, respectivement, qui v\'erifient la
condition {\rm (D)} de {\rm  (2.2)},
alors il existe une unique extension de $\tau  $ \`a une
transformation naturelle
$$
\tau  ' :F'\longrightarrow G'.
$$

Si en outre $\tau  $ est un isomorphisme de foncteurs,
l'extension $\tau  '$ est aussi un isomorphisme. En particulier,
l'extension $G'$ du th\'eor\`eme {\rm (2.2)} est essentiellement unique.
 \endproclaim

\demo{Preuve} La d\'emonstration se fait aussi par r\'ecurrence,
et dans ce cas elle est imm\'ediate d'apr\`es (S4) et
(D).
\enddemo

\noindent$\bold{ (2.5) Variantes.}$
Dans les applications du th\'eor\`eme que nous donnerons dans les
paragraphes suivants, les cat\'egories $\bold{M} $ et
$\bold{W}$ seront diff\'erentes selon les cas, mais on
disposera
toujours
d'un th\'eor\`eme de r\'esolution, d'un lemme de Chow, et la th\'eorie
des hyperr\'esolutions cubiques
sera applicable (cf. \cite{HC}(I.3.11)). Dans toutes ces situations,
la preuve est alors  analogue \`a celle donn\'ee dans le cas
consider\'e ci dessus.

Notons aussi qu'il y a une  version homologique des
th\'eor\`emes (2.2) et (2.4), qui s'obtient par  passage \`a la
cat\'egorie oppos\'ee de la cat\'egorie de
descente. Bien s\^ur, la version homologique, qui
consid\`ere
uniquement des foncteurs covariants, est plus naturelle dans un cadre
abstrait, mais nous avons pr\'ef\'er\'e donner plut\^ot l'enonc\'e
cohomologique
pour des raisons historiques.

\subheading{3. Application \`a l'homotopie de De Rham alg\'ebrique}

La premi\`ere application que nous consid\'erons est \`a
l'homotopie de De Rham des vari\'et\'es alg\'ebriques,
lisses ou non, sur un corps de caract\'eristique z\'ero.

\noindent $\bold{(3.1)}$
Soient $\bold{Top}_{*}$ la cat\'egorie des espaces topologiques
point\'es
et  $Ho\,\bold{Adgc}_{*}(\Bbb{Q})$ la cat\'egorie des
$\Bbb{Q}$-alg\`ebres dgc augment\'ees, localis\'ee
par
rapport aux homologismes. Rappelons que Sullivan (\cite{Su}) a prouv\'e
l'existence d'un foncteur $$
A_{Su}:\bold{Top}_{*}\longrightarrow Ho\,
\bold{Adgc}_{*}(\Bbb{Q}) $$
qui associe \`a tout espace topologique point\'e $(X,x)$ une alg\`ebre de
formes
diff\'erentielles
dont la cohomologie est
isomorphe \`a la cohomologie rationnelle singuli\`ere de $X$, et
qui contient aussi des informations sur le type d'homotopie
rationnelle de $X$, si $X$ est de type fini.
En effet, l'alg\`ebre de Lie $\pi
_{1}(A_{Su}(X,x))$, duale de l'espace des
ind\'ecomposables de degr\'e 1 d'un mod\`ele minimal de $A_{Su}(X,x)$,
est isomorphe \`a l'alg\`ebre de Lie
rationnelle associ\'ee au groupe fondamental de $(X,x)$, c'est-\`a-dire,
l'alg\`ebre
de Lie du
complet\'e de Malcev $\pi _{1}^{Mal}(X,x)$ de $\pi _{1}(X,x)$.
Et, si $X$ est
simplement connexe  et $n\ge 2$, l'espace
$\pi_{n}(A_{Su}(X))$, dual de l'espace d'ind\'ecomposables
de degr\'e $n$ d'un mod\`ele minimal de $A_{Su}(X)$, est isomorphe
au n-i\`eme espace d'homotopie
rationnelle de $X$, $\pi _{n}(X)\bigotimes_{\Bbb{Z}}\Bbb{Q}$.

Si $X$ est une vari\'et\'e alg\'ebrique affine
 non singuli\`ere sur un sous-corps $k$ de $\Bbb{C}$, on peut
r\'ealiser alg\'ebriquement  la construction d'une $k$-alg\`ebre dgc de
Sullivan, car les formes diff\'erentielles alg\'ebriques d\'efinissent
une $k$-alg\`ebre dgc $\Omega^{*}(X)$ telle que
$\Omega^{*}(X)\bigotimes_{k}\Bbb{C}$ est quasi-isomorphe
\`a $A_{Su}\left( X^{an}\right)\bigotimes
_{\Bbb{Q}}\Bbb{C}$, d'apr\`es le th\'eor\`eme de Grothendieck (\cite{G},
voir aussi \cite{Ha} ).

Nous allons g\'en\'eraliser cette r\'ealisation \`a toutes les vari\'et\'es
alg\'ebriques s\'epar\'ees et de type fini sur $k$ (cf. \cite{HC}(III.1.7)).

\proclaim{(3.2) \, Th\'eor\`eme} Soient $k$ un corps de caract\'eristique
zero, et  $Ho \,\bold{Adgc}(k)$  la cat\'egorie des $k$-alg\`ebres dgc,
localis\'ee par rapport aux homologismes. Alors il existe un foncteur
contravariant
$$
A_{DR}:\bold{Sch}(k)\longrightarrow Ho \, \bold{Adgc}(k)
$$
tel que:
\roster
\item Si $X$ est un $k$-sch\'ema  affin et lisse, $A_{DR}(X)$ est
l'alg\`ebre $\Omega^{*} (X)$ des formes diff\'erentielles sur $X$.
\item $A_{DR}$ v\'erifie la propri\'et\'e de descente {\rm(D)}.
\item La cohomologie $H^{*}A_{DR}(X)$ est
isomorphe \`a la cohomologie de De Rham alg\'ebrique de $X$
$
H^{*}_{DR}(X,k)$.
\endroster
En outre, ce foncteur est detemin\'e, \`a quasi-isomorphisme pr\`es, par les
conditions (1) et (2).
 \endproclaim

\demo{Preuve} Ce r\'esultat est une cons\'equence imm\'ediate du th\'eor\`eme
(2.2).
En effet, on prend comme cat\'egorie de descente la cat\'egorie
$\bold{Adgc}(k)$ avec la structure de
cat\'egorie de descente d\'efinie par les donn\'ees suivantes:
La classe
${E}$ est celle des homologismes, et le foncteur simple cubique est
obtenu en utilisant le foncteur simple de Thom-Whitney  $\bold{s}_{TW}$
d\'efinit dans \cite{N}(3.2), compte tenu de la remarque (1.5).
Les conditions
de cat\'egorie de descente sont une cons\'equence du lemme (1.10),
appliqu\'e au
foncteur d'oubli de la structure multiplicative,
car la condition (1.11.2)
est v\'erifi\'ee par le morphisme de comparaison $I$ de
\cite{N}(3.3).

Maintenant, le foncteur $G$ est d\'efini,
pour tout k-sch\'ema
 X
de $\bold{Sch_{Reg}}(k)$, par
$G(X):= \Bbb{R}_{TW}\Gamma (X, \Omega ^{*}_{X})$, o\`u $\Bbb{R}_{TW}$ est
le foncteur d\'eriv\'e dans le sens de $k$-alg\`ebres dgc (\cite{N}(4.4)).
 Les conditions (F1) et (F2) sont ais\'ement  v\'erifi\'ees. D'apr\`es
le th\'eor\`eme de comparaison \cite{N}(3.3), on a un quasi-isomorphisme de
complexes $G(X)\cong DR^{*}(X) $,
o\`u $DR^{*}(X)$ est le
complexe de De Rham ordinaire (\cite{G} ou \cite{Ha}),
d\'efini par
$$
DR^{*}(X)=\Bbb R \Gamma (X,\Omega^{*}_{X}).
$$
Comme la condition  (F3) est une condition
cohomologique,
 pour la v\'erifier
 il suffit de le faire  pour le
 foncteur $DR^{*}$. Dans ce cas (F3) est une cons\'equence
de \cite {Gr}(VI,1.2.1), r\'esultat que nous rappelons dans le lemme
(3.3) ci-dessous.
Alors on peut appliquer le th\'eor\`eme (2.2), qui entra{\^\i}ne
les propri\'et\'es (1) et (2)
de (3.2). Finalement, la propri\'et\'e (3) r\'esulte aussit\^ot
du th\'eor\`eme de comparaison, de (2.4) et de
la suite
exacte en cohomologie de De Rham d'un morphisme birationnel
(\cite{Ha}(II.4.4)).
\enddemo

\proclaim{(3.3) Lemme}
Soit $X_{\bullet}$ le diagramme cart\'esien de $k$-sch\'emas lisses
$$
\CD
\widetilde Y @>  j >> \widetilde X \\
@V g VV @V f VV\\
Y @> i >> X
\endCD $$
o\`u
$j$ et $i$ sont des immersions ferm\'ees, et $f$
est  l'\'eclatement du centre $Y$. Alors,
pour tout $p \ge 0$, le morphisme
$$
\Omega _{X}^{p} \longrightarrow \bold{s}\left(
\Bbb{R}i_{*}
\Omega _{Y}^{p}\bigoplus
\Bbb{R}f_{*}\Omega _{\widetilde X}^{p}
\longrightarrow \Bbb{R}(f\circ j) _{*}\Omega _{\widetilde Y
}^{p}  \right)
$$
est un quasi-isomorphisme.
\endproclaim
\demo{Preuve } Le probl\`eme \'etant local en $X$ pour la topologie \'etale,
on peut supposer que $X$
est l'espace affin $\Bbb{A}_{k}^{m}\times\Bbb{A}_{k}^{n}$ et $Y$ est
$\Bbb{A}^{m}_{k}\times 0$. Alors on applique \cite{Gr}(VI,1.2.1).
\enddemo

\noindent $\bold{(3.4)}$ Soit $x$ un $k$-point d'un $k$-sch\'ema $X$
g\'eom\'etriquement connexe.
Par fonctorialit\'e, le morphisme d'inclusion $\{x\}\longrightarrow X$
induit une augmentation $A_{DR}(X)\longrightarrow k$, dans la
cat\'egorie
homotopique. Notons $A_{DR}(X,x)$ cette alg\`ebre augment\'ee (au
sens homotopique), d'apr\`es la th\'eorie de Sullivan (\cite{Su})
on peut associer \`a $(X,x)$ un mod\`ele minimal de $A_{DR}(X,x)$, ses
espaces d'ind\'ecomposables et leurs duals, $\pi_{n}(A_{DR}(X,x))$,
pour tout $n\ge
1$, qui sont ind\'ependants du mod\`ele minimal et fonctoriels en
$(X,x)$.
On pose $\pi _{n}(X,x)_{DR}=\pi_{n}(A_{DR}(X,x)), \, n\ge
1$. L'espace
$\pi _{1}(X,x)_{DR}$ est muni naturellement d'une structure d'alg\`ebre
de Lie sur $k$ qui est pro-nilpotente, ou, ce qui est \'equivalent,
 d'une structure de sch\'ema en groupes
sur $k$, pro-alg\'ebrique unipotente (voir \cite{D2}).

Le corollaire suivant est alors une cons\'equence de
(3.2), (2.4) et le th\'eor\`eme de comparaison de Grothendieck (\cite{G}).

 \proclaim{Corollaire}
 Soient $k$ un
sous-corps de $\Bbb{C}$,  $X$ un $k$-sch\'ema et $x$ un $k$-point de $X$.
Pour toute immersion  $\sigma :k\rightarrow \Bbb {C}$,
notons $X_{\sigma }^{an}$  l'espace analytique associ\'e \`a
$X\bigotimes_{k,\sigma }\Bbb {C}$. Alors:
\roster
 \item
On a un isomorphisme naturel
$A_{DR}(X)\bigotimes_{k,\sigma }\Bbb{C}\cong
A_{Su}(X_{\sigma }^{an})\bigotimes_{\Bbb{Q}}\Bbb{C}$.
\item On a un isomorphisme naturel
$$\pi _{1}(X,x)_{DR}\otimes_{k,\sigma
}\Bbb{C} \cong \pi _{1}^{Mal}(X_{\sigma
}^{an},x)\otimes_{\Bbb{Q}}\Bbb{C}.
$$
\item Si $X_{\sigma }^{an}$ est
simplement connexe, on a des isomorphismes naturels
$$
\pi _{n}(X,x)_{DR}\otimes_{k,\sigma }\Bbb{C}\cong \pi
_{n}(X_{\sigma }^{an})\otimes_{\Bbb{Z}}\Bbb {C}    ,
$$
pour tout entier $n\ge 2$.
\endroster
\endproclaim

\subheading{4. Application au complexe filtr\'e de Hodge-De Rham}

\noindent$\bold{(4.1)}$
Soit $X$  une vari\'et\'e analytique complexe, et soit
$\Omega ^{*}_{X}$ le complexe de De Rham des faisceaux des formes
diff\'erentielles holomorphes sur $X$. Il est bien connu que ce complexe
est muni de la filtration de Hodge $F^{p}\Omega
^{*}(X)=\Omega ^{*\ge p}(X)$, et qu'il en r\'esulte la suite spectrale de
Hodge-De Rham
$$
\oplus_{p+q=r}^{}H^{q}(X,\Omega ^{p}_{X}) \Longrightarrow
 H^{r}(X,\Omega ^{*}_{X})\cong H^{r}(X,\Bbb{C}) , $$
laquelle, par le th\'eor\`eme de Hodge, d\'eg\'en\`ere au terme $E_{1}$ si $X$
est une vari\'et\'e k\"ahlerienne
compacte.

 Dans le cadre alg\'ebrique, Du Bois a prouv\'e (\cite{DB}, voir aussi
\cite{HC}(V.3.5)), en utilisant la th\'eorie de Hodge-Deligne, que si $X$
est une vari\'et\'e alg\'ebrique
sur $\Bbb{C}$, il existe un complexe filtr\'e de Hodge-De Rham
$(\underline{\Omega}^{*} _{X},F)$  qui est une g\'en\'eralisation naturelle
au cas possiblement singulier du complexe pr\'ec\'edent, en
particulier on a une suite spectrale
$$
\oplus_{p+q=r}^{}H^{q}(X,\text{Gr}_{F}^{p}\underline{\Omega}
^{*}_{X}[p]) \Longrightarrow
 H^{r}(X,\underline{\Omega} ^{*}_{X})\cong H^{r}(X,\Bbb{C}) , $$
laquelle, par la th\'eorie de Hodge-Deligne, d\'eg\'en\`ere au terme $E_{1}$ si
$X$  est une vari\'et\'e alg\'ebrique compacte.

Nous prouvons dans ce paragraphe, comme application du th\'eor\`eme (2.2)
et sans recours \`a la th\'eorie de Hodge-Deligne, l'existence du
complexe filtr\'e de Hodge-De Rham $(\underline{\Omega}_{X}^{*},F)
 $ pour tout espace analytique (toujours r\'eduit, s\'epar\'e et
d\'enombrable \`a l'infini).

\proclaim{(4.2) Th\'eor\`eme}
Pour tout espace analytique complexe $X$, il existe un complexe
filtr\'e de faisceaux de $\Bbb{C}$-espaces vectoriels
$(\underline{\Omega }^{*}_{X},F)$,
dont les diff\'erentielles  sont des op\'erateurs diff\'erentiels d'ordre
$\le \, 1$, dont les gradu\'es sont des complexes de
$\Cal{O}_{X}$-modules \`a cohomologie coh\'erente,
et qui v\'erifie les propri\'et\'es suivantes:
\roster
\item Soit $(\Omega _{X}^{*},\sigma )$ le complexe des diff\'erentielles
de K\"ahler de $X$, filtr\'e par la filtration b\^ ete $\sigma $, alors il
existe un morphisme naturel de complexes filtr\'es
$$
(\Omega _{X}^{*},\sigma )\longrightarrow (\underline{\Omega
}_{X}^{*},F),
$$
qui est un quasi-isomorphisme filtr\'e si $X$ est une vari\'et\'e complexe.
 \item Si $f:X\longrightarrow Y$ est un morphisme
d'espaces analytiques, il existe un morphisme naturel
$$
f^{*}:(\underline{\Omega}_{Y}^{*},F)\longrightarrow
\Bbb{R}f_{*}(\underline{\Omega}_{X}^{*},F).
$$
\item Le complexe   $(\underline{\Omega}_{X}^{*},F)$ v\'erifie la
propri\'et\'e de descente {\rm(D)}.
\item
Le complexe $\underline{\Omega}^{*}_{X}$  est une r\'esolution de
$\Bbb{C}_{X}$.
\item Si $X$ est un $\Bbb{C}$-sch\'ema, le complexe
$(\underline{\Omega}_{X^{an}}^{*},F)$ est naturellement isomorphe au
complexe d\'efini par Du Bois.
\item  $Gr^{p}_{F}\underline{\Omega }^{*}_{X}=0$, si $p\notin
[0,\text{dim }X]$.
 \endroster
En outre, ce complexe est d\'etermin\'e, \`a
quasi-isomorphisme filtr\'e pr\`es, par les conditions {\rm (1)} \`a {\rm
 (3)}.
 \endproclaim

\demo{Preuve} Nous utiliserons une variante en g\'eometrie analytique  du
th\'eor\`eme (2.2). Soient
$\bold{W}$ la cat\'egorie des espaces analytiques propres sur $X$,
$\bold{M}$  la sous-cat\'egorie pleine
de $\bold{V}$  des vari\'et\'es complexes,
 $\bold{D}$ la cat\'egorie
des complexes filtr\'es de faisceaux de $\Bbb{C}$-espaces vectoriels
sur $X$ dont les
gradu\'es sont des complexes de $\Cal{O}_{X}$-modules  \`a
cohomologie coh\'erente,
et leurs  diff\'erentielles sont des operateurs diff\'erentiels
d'ordre $\le\,
1$. Nous munissons $\bold{D}$  de la structure de cat\'egorie de descente
suivante:
La classe ${E}$ est form\'ee par les  quasi-isomorphismes
filtr\'es et le foncteur simple de complexes se filtre par la filtration
$(s,s)$ de \cite{N}(6.1),
c'est-\`a-dire, si $(C,F)$ est un objet cubique de
$\bold{D}$, on d\'efinit $
(\bold{s,s})(C,F)=(\bold{s}C,\bold{s}F).
$
 Les conditions de cat\'egorie de  descente sont
une cons\'equence du lemme
(1.10), le foncteur $\phi $ \'etant le foncteur de passage au
gradu\'e $Gr= \prod_{n} Gr_{n} $.

Soit $G(Z)\,=\,\Bbb{R}f_{*}(\Omega^{*}_{Z},F)$,   o\`u $F$
d\'enote
la filtration de Hodge, pour tout objet  $f:Z\longrightarrow X$ de
$\bold{M}$.
Alors $G$ est un foncteur contravariant
\`a valeurs dans $\bold{D}$.

En outre $G$ est un
foncteur v\'erifiant (F1) \`a (F3). En effet, (F1) et (F2) sont
triviales,
et la propri\'et\'e (F3) se suit de \cite{Gr}(VI,1.2.1), r\'esultat que
nous avons rappel\'e dans le lemme (3.3), car la filtration
de Hodge dans le cas non
singulier est la filtration par le degr\'e, et  on a donc
$Gr^{p}_{F}\Omega ^{*}=\Omega
^{p}[-p]$.

Il est claire qu'avec ces donn\'ees les conclusions des th\'eor\`emes
(2.2) et (2.4) sont encore  v\'erifi\'ees. En effet,
la preuve est toute \`a fait analogue au cas alg\'ebrique, car
d'apr\`es le th\'eor\`eme de r\'esolution de singularit\'es dans le cas
analytique
(\cite{AH} , voir  \cite {BM})
on a aussi dans ce contexte les r\'esultats sur les hyperr\'esolutions
cubiques n\'ecessaires (\cite{GN}(I.3.11.2)).
 En outre, dans le point (2.3.1) de la preuve,  on
doit
prouver que si $X_{\bullet}\longrightarrow X$ est une hyperr\'esolution
cubique d'une vari\'et\'e complexe $X$, $G(X)\longrightarrow G(X_{\bullet})$
est un quasi-isomorphisme filtr\'e, mais cette condition est local dans
$X$ et on applique le lemme de Chow local de Hironaka (\cite{H1}).

La propri\'et\'e (2) est une cons\'equence de l'existence du complexe de
Hodge-De Rham dans le cas du diagramme d\'efini par un morphisme, compte
tenu de la commutativit\'e avec le foncteur {\rm tot }.
Le m\^eme argument donne aussi la naturalit\'e. En fait, on
peut aussi  obtenir
la fonctorialit\'e avec une version du th\'eor\`eme
d'extension pour une
cat\'egorie fibr\'ee en cat\'egories de descente au-dessus d'une
cat\'egorie, mais que nous ne formulons pas.

 La propri\'et\'e (4) r\'esulte de
(2.4) et \cite{GN}(I.6.9), (5) se d\'emontre comme dans \cite{GN} (V.3.7),
et  (6) r\'esulte de l'existence d'une hyperr\'esolution
$X_{\bullet}\longrightarrow X$  de $X$ telle que $\text{dim }X_{\alpha
}\le \text{dim }X - \vert \alpha \vert +1$ (\cite{GN}(I.2.15)).
\enddemo

\noindent $\bold{(4.3)}$
Avec ce complexe filtr\'e de Hodge-De Rham $(\underline{\Omega
}_{X}^{*},F)$ on peut \'etendre la th\'eorie classique de Hodge-De Rham aux
espaces analytiques de la fa\c con suivante:

D'abord, pour tout espace analytique $X$, on d\'efinit les complexes de
faisceaux
$$
\underline{\Omega
}^{[p]}_{X}:=\text{Gr}^{p}_{F}\underline{\Omega}_{X}^{*}[p] ,\, 0\le
p\le \text{dim }X,
$$
 qui sont \`a cohomologie coh\'erente et fonctoriels en
$X$.

Ensuite, on d\'efinit la cohomologie de Hodge de $X$ par
$$
\oplus_{}^{}H^{p,q}(X):=\bigoplus_{}^{}\Bbb{H}^{q}(X,\underline{\Omega
}^{[p]}_{X}), \, 0\le p,q \le n,
$$
qui est fonctorielle en $X$ et, si $X$ est compact, est de dimension
finie, donc on peut d\'efinir dans ce cas les nombres de Hodge de $X$ par
$$
h^{p,q}(X)= \text{dim}_{\Bbb{C}}H^{p,q}(X).
$$

Finalement on a une suite spectrale de Hodge-De Rham
associ\'ee au complexe filtr\'e $(\underline{\Omega}_{X}^{*},F)$,
$$
\oplus_{p+q=r}^{}H^{p,q}(X) \Longrightarrow H^{r}(X,\Bbb{C}).
$$
d'o\`u, si $X$ est compact,
$$
\chi (X)=\sum_{p,q}^{}(-1)^{p+q}h^{p,q}(X)
$$
est la caract\'eristique d'Euler de $X$. Et, si $X$ est compact
et sous-k\"ahlerien (i.e. il existe une
vari\'et\'e k\"ahlerienne $X'$ et un morphisme surjectif propre
$X'\longrightarrow X$), par la th\'eorie de Hodge-Deligne, cette suite
spectrale d\'eg\'en\`ere au terme $E_{1}$.

\subheading{5. Application \`a la th\'eorie des motifs}

\noindent $\bold{(5.1)}$
Soit $k$ un corps de caract\'eristique z\'ero. Nous noterons
$\bold{V}(k)$ la cat\'egorie des $k$-sch\'emas projectifs et
lisses, et
$\Cal{M}^{+}_{rat}(k)$ la cat\'egorie des motifs de Chow effectifs sur
$k$ (voir \cite{M},
\cite{Sc}). On rappelle que cette derni\`ere cat\'egorie est
pseudo-ab\'elienne,
et qu'on a un foncteur contravariant
$$
h:\bold{V}(k)\longrightarrow
\Cal{M}^{+}_{rat}(k),
$$
qui associe \`a tout $k$-sch\'ema projectif et
lisse $X$, le motif effectif $h(X)=(X,id_{X})$.
L'objectif de ce paragraphe est \'etendre le
foncteur $h$  a tous les $k$-sch\'emas.

On peut faire cette extension de deux
fa\c cons diff\'erents. La plus imm\'ediate c'est comme une th\'eorie \`a
support
compact, mais, comme on verra, on peut faire aussi l'extension \`a une
th\'eorie sans supports.
D'abord nous prouvons ci-dessous l'existence d'une
extension du foncteur
$h$ aux  sch\'emas propres sur $k$.

\noindent$\bold{(5.2)}$
Notons  $\bold{Sch_{Prop}}(k)$ la cat\'egorie des sch\'emas propres sur $k$,
et $C^{b}(\Cal{M}^{+}_{rat}(k))
$ la
cat\'egorie des complexes born\'es de motifs de Chow effectifs sur $k$.

\proclaim{Proposition}
Il existe une unique extension de $h$ \`a un foncteur
contravariant $$
(\Phi^{op} ,\bold{Sch_{Prop}}(k))\longrightarrow  Ho \,
(\Phi , C^{b}(\Cal{M}^{+}_{rat}(k))),
$$
not\'e encore $h$,
qui commute avec le foncteur $tot $, et qui
 v\'erifie la propri\'et\'e de descente {\rm (D)}.
\endproclaim

\demo{Preuve} Nous utiliserons une variante du th\'eor\`eme (2.2)
obtenu dans la situation suivante.

 Soient $\bold{M}=\bold{V}(k)$,
$\bold{W}=\bold{Sch_{Prop}}(k)$, et
munissons $C^{b}(\Cal{M}^{+}_{rat}(k))$
de la structure de
cat\'egorie de descente d\'efinie dans (1.8).
Le foncteur $h:\bold{V}(k)\longrightarrow
\Cal{M}^{+}_{rat}(k)
$ induit un foncteur contravariant
$$h:(\Phi ^{op},\bold{V}(k))\longrightarrow
Ho \left(\Phi  , C^{b}(\Cal{M}^{+}_{rat}(k))\right)
$$
qui v\'erifie
trivialement les conditions (F1) et (F2) de (2.2).
D'apr\`es \cite{M},
$h$ v\'erifie aussi la propri\'et\'e (F3).
En effet, ceci d\'ecoule du lemme (5.3) qui sera prouv\'e ci-dessous.

D'apr\`es le lemme de Chow et le
th\'eor\`eme de r\'esolution, (2.2) est v\'erifi\'e aussi avec ces hypoth\`eses
(voir
\cite{HC}(3.11.3)),donc  $h$ s'\'etend de fa\c con unique \`a un foncteur
$$(\Phi ^{op},\bold{W})\longrightarrow Ho\,\left(
{\Phi },C^{b}(\Cal{M}^{+}_{rat}(k))\right) $$
qui v\'erifie la propri\'et\'e de descente (D).

\enddemo

\proclaim{(5.3) \, Lemme}
Soient $i:Y\longrightarrow X$ une immersion ferm\'ee de $k-$sch\'emas
projectifs et lisses, $f:\widetilde X\longrightarrow X$ l'\'eclatement
de $X$ le long de $Y$, et
$$
 \CD
\widetilde Y @>j>> \widetilde X \\
@V g VV @VV f V\\
Y @>i>> X
\endCD
$$
le diagramme cart\'esien construit \`a partir de $f$ et $i$.
Alors la suite de motifs de Chow
 $$
0\longrightarrow h(X)@> i^{*}+f^{*} >>
h(Y)\bigoplus
h(\widetilde X)@> g^{*}-j^{*}>> h(\widetilde
Y)\longrightarrow 0 $$
est exacte et scind\'ee, en particulier le morphisme
$$
h(X)\longrightarrow \bold{s}\left(
h(Y)\bigoplus
h(\widetilde X)\longrightarrow  h(\widetilde
Y)\right)
$$
est un homotopisme.
\endproclaim

\demo{Preuve} Ceci
est une cons\'equence du calcul effectu\'e par Manin du motif de Chow d'un
\'eclatement
(\cite{M}, \S 9 Cor., voir aussi \cite{Sc}, Th.(2.8)). En effet,
d'apr\`es loc.cit. on a un isomorphisme $$
\varphi :h(X)\oplus_{}^{}\left(\bigoplus_{i=1}^{r-1} h(Y)(-i)
\right)\longrightarrow  h(\widetilde X),
$$
d\'efini par
$$
\varphi=f^{*}+\Sigma
_{i=1}^{r-1}j_{*}\circ(\xi ^{i-1}\cup g^{*}), $$
o\`u $\xi $  est
la premi\`ere classe de Chern de $\Cal{ O}_{\widetilde Y}(1). $
On a aussi l'isomorphisme
$$
\psi:\bigoplus_{i=0}^{r-1}h(Y) (-i)\longrightarrow h(\widetilde Y)
$$
d\'efini par
$$
\psi = \Sigma
_{i=0}^{r-1}\xi ^{i}\cup g^{*}
$$
Donc, avec ces isomorphismes, si on pose
$$
A=\bigoplus_{i=1}^{r-1} h(Y)(-i),
$$
 la suite du lemme s'\'ecrit
simplement
 $$
0\longrightarrow h(X)@> i^{*}\times id_{X}\times 0
 >>
h(Y)\oplus
h(X)\oplus_{}^{} A
@>(id_{Y}-i^{*}+0)  \times (0+ 0 -id_{A}) >>
h(Y)\oplus A
\longrightarrow 0
$$
laquelle est trivialement exacte et scind\'ee.
\enddemo

\noindent $\bold{(5.4)}$
Maintenant nous donnons  l'extension du foncteur $h$
correpondante \`a une th\'eorie \`a support compact.
Pour ceci notons  $\bold{Sch}_{c}(k)$ la cat\'egorie des sch\'emas
sur $k$, avec des morphismes propres.

\proclaim{ \, Th\'eor\`eme}
Il existe un unique foncteur
contravariant $$
h_{c}:\bold{Sch}_{c}(k)\longrightarrow  Ho \,
C^{b}(\Cal{M}^{+}_{rat}(k)) $$
tel que:
\roster
\item Si $X$ est projectif et lisse sur $k$, $h_{c}(X)$ est le motif de
Chow associ\'e \`a $X$.
\item $h_{c}$ v\'erifie la propri\'et\'e de descente {\rm (D)}.
\item Si $Y$ est un sous-sch\'ema ferm\'e d'un sch\'ema  $X$, on a un
isomorphisme $h_{c}(X - Y)\cong\bold{s}(h_{c}(X) \rightarrow
h_{c}(Y))$ \endroster
\endproclaim

\demo{Preuve} Nous allons d\'efinir $h_{c}$ plus g\'en\'eralement pour tout
objet de $(\Phi ,\bold{Sch}(k))$, ceci donnera comme
cons\'equence la
fonctorialit\'e. Soit $\bold{I}$ une cat\'egorie ordonnable finie, $U_{\bullet
}:\bold{I}^{op}\longrightarrow \bold{Sch}_{c}$ un foncteur.
D'apr\`es \cite{HC}(I.4.2) et (I.4.3), il existe une compactification
$U_{\bullet}\longrightarrow X_{\bullet}$ de $U_{\bullet}$, telle que le
compl\'ementaire $Y_{\bullet}$ est un aussi un diagramme de sch\'emas.
On d\'efinit
 $h_{c}(U_{\bullet}):=\bold{s}\left(h(Y_{\bullet}\rightarrow
X_{\bullet})\right)$.
D'apr\`es (D), et \cite{HC}(I.4.4), $h_{c}(U_{\bullet})$
ne d\'epend pas de la compactification et on obtient un foncteur $h_{c}$
qui v\'erifie les propri\'et\'es (1) et  (2) du
th\'eor\`eme. La propri\'et\'e (3) r\'esulte ais\'ement de la
d\'efinition de $h_{c}$ et de la propri\'et\'e de descente (D). L'unicit\'e
est imm\'ediate.
\enddemo

 \noindent$\bold{(5.5)}$
Nous rappelons  que dans une cat\'egorie pseudo-ab\'elienne tout complexe
born\'e contractile  $L^{*}$
a une caract\'eristique d'Euler $\chi (L^{*})=\Sigma (-1)^{i}[L^{i}]$
nulle dans le groupe  ${K}_{0}$ correspondant.
Pour la commodit\'e du lecteur nous donnons ici une preuve de ce
r\'esultat \`a partir de la proposition suivante.

\proclaim{Proposition} Soit $\bold{A}$ une cat\'egorie
pseudo-ab\'elienne. Si $L$ est
un complexe contractile de $C(\bold{A})$, il existe un complexe \`a
diff\'erentielle nulle $P$
et un isomorphisme de complexes  $\varphi :Con(P)\longrightarrow L$, o\`u
$Con(P)$ est le complexe c\^one de $P$. \endproclaim
\demo{Preuve}
Soit $h$ une contraction de $L$, c'est-\`a-dire,
$h:L\longrightarrow L[-1]$ est un morphisme gradu\'e tel que $1=hd+dh$.
Alors $dh$ et $hd$ sont des projecteurs compl\'ementaires, et puisque
$\bold{A}$ est pseudo-ab\'elienne,
$L$ est isomorphe, en tant qu'objet gradu\'e, \`a la somme $Ker \, dh
\,\oplus\, Im \, dh$, avec les
\'egalit\'es $Ker \, dh \,=\, Im \, hd$ et $Im \, dh\, =\, Ker \, hd$.
Puisque $d^{2}=0$, la restriction de $d$ \`a $Im \, dh$ est nulle,
et la
restriction de $d$ \`a $Ker \, dh$ induit un isomorphisme de complexes
$Ker \, dh\longrightarrow Im \, dh [+1]$. Alors $P= Im \, dh$
v\'erifie les conditions du lemme.
\enddemo

\proclaim{(5.6) \, Corollaire} Soit $\bold{A}$ une cat\'egorie
pseudo-ab\'elienne, $C^{b}(\bold{A})$  la cat\'egorie de complexes born\'es de
$\bold{A}$
et  $Ho \,C^{b}(\bold{A})$ sa localisation par rapport aux
homotopismes. Alors la caract\'eristique d'Euler
$$
\chi :Ob \,  C^{b}(\bold{A})\longrightarrow {K}_{0}(\bold{A})
$$
induit une application sur la localisation
$$
\chi :Ob\,  Ho \,  C^{b}(\bold{A})\longrightarrow
{K}_{0}(\bold{A}). $$
\endproclaim

\noindent
$\bold{(5.7)}$ Puisque la cat\'egorie des motifs de Chow est
pseudo-ab\'elienne,
on d\'eduit de (5.4) et (5.6) le
r\'esultat suivant, qui r\'esout le probl\`eme de Serre cit\'e, et qui
avait \'et\'e d\'ej\`a prouv\'e par Gillet et Soul\'e (\cite{GS}).

\proclaim{Corollaire} Il existe une application unique
$$
\chi _{c}:Ob \,  \bold{Sch}(k) \longrightarrow
{K}_{0}(\Cal{M}^{+}_{rat}(k)) $$
qui v\'erifie:
\roster
\item $\chi _{c}(X)=[X]$, si $X$ est un sch\'ema projectif et lisse.
\item $\chi _{c}(X - Y)= \chi _{c}(X)-\chi _{c}(Y)$, si $Y$ est un
sous-sch\'ema  ferm\'e de $X$.
\endroster
\endproclaim
\noindent $\bold{(5.8)}$
Maintenant nous
consid\'erons l'extension de $h$ qui correspond \`a une th\'eorie sans
supports. Rappelons pour ceci  que si $\Cal{M}_{rat}(k)$ est la
cat\'egorie des motifs de Chow, on a aussi
un foncteur covariant $h_{*}:\bold{V}(k)\longrightarrow
\Cal{M}_{rat}(k)$
tel que $h_{*}(X)=h(X)(dim \,X)$, si $X$ est une vari\'et\'e
projective et lisse. Si $f:X\longrightarrow Y$
est un morphisme entre telles vari\'et\'es,
on appelle $f_{*}:h_{*}(X)\longrightarrow h_{*}(Y)$ le
morphisme de Gysin induit par $f$.
Pour cette th\'eorie covariante nous consid\'ererons aussi la categorie
$\Cal{ C}^{b}(\Cal{M}_{rat}(k))$ comme une cat\'egorie de descente
homologique, avec le foncteur simple ordinaire des complexes de cha{\^\i}nes,
qui nous noterons encore $\bold{s}$.

Si $X$  est un sch\'ema  projectif et lisse sur $k$  et
 $Y$ est un diviseur \`a croisements normaux dans $X$, qui est une
r\'eunion de diviseurs lisses $Y=\bigcup_{\alpha =1}^{r}Y_{\alpha }$,
on a un diagramme cubique augment\'e  $S_{\bullet}(Y)\longrightarrow X$,
d\'efini par les intersections des $Y_ {\alpha }$, $S_{\chi
}(Y):=\bigcap_{\chi  (\alpha )=1}
Y_{\alpha }$, pour $\chi \in \square_{r}$. Par la fonctorialit\'e de
$h_{*}$ on
obtient un objet cubique
$h_{*}(S_{\bullet}(Y)\longrightarrow X)(-dim \, X)$ de
$\Cal{M}^{+}_{rat}(k)$
et donc un  complexe
$\bold{s}(h_{*}(S_{\bullet}(Y)\longrightarrow X))(-dim \, X)$ dans
$\Cal{C}^{b}(\Cal{M}^{+}_{rat}(k))$, que nous noterons
$G_{\bullet}(X,Y)$. Ce complexe de Gysin du couple $(X,Y)$ est
l'analogue dans le pr\'esent contexte
du terme $E_{1}$ de la suite spectrale d\'eduite du complexe de De Rham
logarithmique filtr\'e par la filtration par le poids (\cite{D1}
(3.2.4), \cite{GN}(\S 1)).

Notons $\bold{P}$ la cat\'egorie des couples $(X,Y)$, o\`u $X$ est un
sch\'ema projectif et lisse sur $k$ et  $Y=\bigcup_{\alpha
=1}^{r}Y_{\alpha }$ est un  diviseur  \`a  croisements normaux dans $X$,
qui est r\'eunion de diviseurs lisses.  Notons $I$ l'ensemble ordonn\'e
$[1,r]$.
 Un morphisme $f:(X',Y')\longrightarrow (X,Y)$ est un morphisme
de sch\'emas $f:X'\longrightarrow X$ tel que, pour tout $\alpha \in I$,
l'image inverse $f^{-1}(Y_{\alpha })$ soit une somme $\sum_{\beta
=1}^{s}
m_{\alpha ,\beta }Y'_{\beta }$ de composantes irr\'eductibles de $Y'$
. Nous noterons $M=(m_{\alpha ,\beta })_{(\alpha ,\beta )\in I\times
J}$ la matrice des multiplicit\'es de $f$. Nous  prouverons dans
(5.11) la
fonctorialit\'e de $G(X,Y)$ (cf. \cite {D1} (8.1.19.2)). Pour ceci nous
utiliserons la variante suivante de la formule d'exc\'es d'intersection
de Fulton-MacPherson.
\proclaim{(5.9)\, Proposition}
Soit
$$
 \CD
 Y' @>j>>  X' \\
@V g VV @VV f V\\
Y @>i>> X
\endCD
$$
 un diagramme cart\'esien de $k$-sch\'emas projectifs tel que $i$
est une immersion ferm\'ee, $X,\,X'\,$ et $Y$ sont lisses et $Y'$ est
un diviseur \`a croisements normaux dans $X'$ dont les composantes
irr\'eductibles $Y'_{\beta }, \, \beta =1,...,s$ sont  lisses.
Notons $\Cal{E}$ le fibr\'e normal d'exc\'es sur $Y'$ d\'efini
par la suite exacte $$
0\longrightarrow \Cal{N}_{Y'/X'}\longrightarrow
g^{*}\Cal{N}_{Y/X}\longrightarrow \Cal{ E} \longrightarrow 0.
$$
Soient
$f^{-1}(Y)=\sum_{\beta =1}^{s}m_{\beta } Y'_{\beta }$
le diviseur image inverse
de $Y$,
$\eta_{\beta }:Y'_{\beta }\longrightarrow Y'$ l'inclusion,
$g_{\beta }:=g\circ \eta_{\beta },\,
j_{\beta }:=j\circ \eta_{\beta }$, $ \Cal{E}_{\beta }
:=\eta_{\beta }^{*}\Cal{E}$ et
$\xi_{\beta }$ la classe de Chern de degr\'e
maximum de  $ \Cal{E}_{\beta }$. Alors le diagramme $$
 \CD
\bigoplus_{\beta =1}^{s} h( Y'_{\beta })(-1) @>\sum_{\beta =1}^{s}
j_{\beta *}>>
h(X') \\ @A \sum_{\beta =1}^{s} m_{\beta }g_{\beta }^{*}\cup
\xi_{\beta } AA @AA f^{*} A\\ h(Y)(dim \, Y-dim \, X) @>i_{*}>>h(X)
\endCD
$$
est commutatif.
\endproclaim

\demo{Preuve} Par le principe d'identit\'e de Manin
(\cite {M}(\S 3), cf. \cite{Sc}(2.3))
il suffit de prouver la commutativit\'e du diagramme pour le groupe de
Chow $A^{*}$. Nous utilisons pour ceci l'anneau de Chow
operationnel d\'efini par
Fulton (\cite{F}(17.3)) pour tout sch\'ema de type fini sur $k$, et qui
co{\"\i}ncide avec l'anneau de Chow classique sur les variet\'es projectives
et lisses.
En particulier, on a la
 "Excess Intersection Formula" (\cite{F}(17.4.1))
$$f^{*}(i_{*}(y))=j_{*}(g^{*}(y) \cup \xi ),$$
pour tout
$y \in A^{*}(Y)$, o\`u $\xi \in A^{e}(Y') $ est la $e$-i\`eme classe de
Chern de $\Cal{E}$.
D'apr\`es \cite{F}(17.4.7)(i), et avec ses notations, on a
$[j]=\sum_{\beta =1}^{s}m_{\beta }\eta_{\beta
*}([j_{\beta }])$, o\`u [\,] denote
la classe d'orientation, et compte tenu de \cite{FM}(I.2.5)(G$_{3}$)(ii)
il r\'esulte $ j_{*}=\sum_{\beta =1}^{s}m_{\beta } j_{\beta *}\circ
\eta^{*}_{j}$, d'o\`u on  d\'eduit
$$
\align
f^{*}(i_{*}(y))=&j_{*}(g^{*}(y )\cup \xi )\\
=& \sum_{}^{}m_{\beta}j_{\beta *}\eta^{*}_{\beta }
(g^{*}(y )\cup \xi)\\
=& \sum_{}^{}j_{\beta *}(m_{\beta } g^{*}_{\beta }(y )\cup
\xi_{\beta } ),
\endalign
$$
pour tout $y \in A^{*}(Y)$, ce qui prouve la
proposition.
\enddemo

\proclaim{(5.10)\, Corollaire}
Soit
$$
 \CD
 Y' @>j>>  X' \\
@V g VV @VV f V\\
Y @>i>> X
\endCD
$$
 un diagramme commutatif de $k$-sch\'emas projectifs tel que $i$
est une immersion ferm\'ee, $X,\,X'\,$ et $Y$ sont lisses, $Y$ est un
diviseur de $X$ et $Y'$
est un diviseur \`a croisements normaux dans $X'$ dont les composantes
irr\'eductibles $Y'_{\beta }, \, \beta =1,...,s$, sont  lisses. Si
$f^{*}\Cal{O}(Y)=\sum_{}^{}m_{\beta }\Cal{O}(Y'_{\beta })$, alors on a
$$
f^{*}\circ i_{*}=\sum_{\beta }^{} m_{\beta } j_{\beta *}\circ
g^{*}_{\beta }. $$
\endproclaim
\demo{Preuve}
Si $Y'$ est $f^{-1}(Y)$ on applique la proposition ant\'erieure. Dans
l'autre cas on a $X'=f^{-1}(Y)$. Il r\'esulte nouvement par le principe
d'identit\'e de Manin et de \cite{F}(17.4.1)
 $$
f^{*}(i_{*}(\alpha ))=
\sum^{}_{\beta }m_{\beta }\xi'_{\beta
}\cup f^{*}(\alpha ),
$$
o\`u $\xi'_{\beta }$ est la classe de Chern de $Y'_{\beta }$. Alors on
applique la  formule de projection et on obtient
$$
f^{*}(i_{*}(\alpha ) )=
\sum_{\beta }^{}m_{\beta }j_{\beta *}g_{\beta }^{*}(\alpha ),
$$
pour tout $\alpha \in A^{*}(Y),$ ce qui prouve le corollaire.
\enddemo

\proclaim{(5.11)\, Proposition} Soit $f:(X',Y')\longrightarrow
(X,Y)$ un morphisme de $\bold{P}$. Il existe un morphisme
de complexes $G_{\bullet}(f):G_{\bullet}(X,Y)\longrightarrow
G_{\bullet}(X',Y')$,
tel que la restriction $G_{0}(f)$ de $G(f)$ \`a $h(X)$ est $f^{*}$.
Ces morphismes v\'erifient $G_{\bullet}(f\circ g)=G_{\bullet}(g)\circ
G_{\bullet}(f)$ et $G_{\bullet}(id)=id. $  \endproclaim
\demo{Preuve}
On pose $G_{0}(f)=f^{*}$. D'abord on d\'efinit
$G_{\bullet}(f)$ sur la composante $G_{1}(X,Y)$
 de la fa\c con
suivante.
Pour toute couple $(\alpha ,\beta )\in I\times J$ telle que
$f(Y'_{\beta })\subseteq Y_{\alpha }$ notons $f_{\alpha ,\beta
}:Y'_{\beta }\longrightarrow Y_{\alpha }$ la restriction de $f$.
Alors on d\'efinit un morphisme
$$G_{\alpha ,\beta }(f):h_{*}(Y_{\alpha })(-dim \, X)\longrightarrow
h_{*}( Y_{\beta }')(-dim\, X') $$
par $G_{\alpha ,\beta }(f)=m_{\alpha ,\beta }f^{*}_{\alpha
,\beta }$.
En composant $G_{\alpha ,\beta }(f)$ avec l'inclusion canonique
$h_{*}(Y'_{\beta })(-dim\, X')\longrightarrow G_{1}(X',Y')$ et
faisant la somme par rapport \`a $\beta $ on obtient un morphisme
$$
 G_{1,\alpha }(f):h_{*}(Y_{\alpha })(-dim\,X)\longrightarrow
G_{1}(X',Y').
$$
D'apr\`es (5.10) on a pour chaque $\alpha  \in I$
$$
G_{0}(f)\circ i_{\alpha *}
=\sum_{\beta =1}^{s}j_{\beta *}\circ
G_{1,\alpha }(f).
 $$
Ces morphismes $G_{1,\alpha }(f)$ sont les composantes d'un morphisme
$$
G_{1}(f):G_{1}(X,Y)\longrightarrow G_{1}(X',Y')
$$
qui v\'erifie
$$
G_{0}(f)\circ \sum_{\alpha =1}^{r}i_{\alpha *}
=\sum_{\beta =1}^{s}j_{\beta *}\circ G_{1}(f)
{}.
$$
Le morphisme $G_{1}(f)$ s'\'etend naturellement \`a toutes les
composantes $G_{p}(X,Y)$.
En effet, pour toute couple de suites croissantes $\sigma=(\sigma
(1),\cdots,\sigma (p)) \subseteq I$
et $\tau =(\tau (1),\cdots, \tau (p))\subseteq J$ telle que
$f(Y'_{\tau })\subseteq Y_{\sigma }$, notons
$f_{\sigma,\tau  }:Y'_{\tau  }\longrightarrow Y_{\sigma  }$ la
restriction de $f$. Alors  on d\'efinit un morphisme $$
G_{\sigma,\tau }(f):h_{*}(Y_{\sigma })(-dim \,
X)\longrightarrow h_{*}(Y'_{\tau })(-dim \, X') $$
par $G_{\sigma,\tau  }(f):=m_{ \sigma,\tau   }f^{*}_{\sigma,\tau  }
$, o\`u $m_{\sigma,\tau }$  est le determinant du mineur d'indices
$(\sigma,\tau  ) $ de la matrice $M$ des multiplicit\'es de $f$. En
composant $G_{\sigma,\tau }(f)$ avec l'inclusion
canonique $ h_{*}(Y'_{\tau })(-dim \, X)\longrightarrow
G_{p}(X',Y') $   et faisant la somme par rapport \`a $\tau $
on obtient les composantes
$$G_{p,\sigma }(f):h_{*}(Y_{\sigma })(-dim\,X)\longrightarrow
G_{p}(X',Y')$$
 du morphisme
$$G_{p}(f):G_{p}(X,Y)\longrightarrow G_{p}(X',Y').$$
La fonctorialit\'e de $G$ sera une cons\'equence de la fonctorialit\'e de
l'alg\`ebre ext\'erieur.
En effet, si $g:(X'',Y'')\longrightarrow (X',Y')$ est un morphisme de
$\bold{P}$, et $M'$ est leur matrice de multiplicit\'es,
alors la matrice $M''$ des multiplicit\'es de la composition $f\circ g$
est le produit des matrices $M'\circ M$. La fonctorialit\'e r\'esulte de
l'identit\'e $m''_{\sigma, \rho }=\sum_{\tau  }^{}m_{\sigma
,\tau  }m'_{\tau ,\rho  }$, pour toute couple $(\sigma  ,\rho  )$, qui
exprime en termes de d\'eterminants la fonctorialit\'e de l'alg\`ebre
ext\'erieur.

Finalement, v\'erifions que les morphismes $G_{p}(f)$ d\'efinissent un
morphisme de complexes
$$ G_{\bullet}(f):G_{\bullet}(X,Y)\longrightarrow G_{\bullet}(X',Y'),
$$
c'est-\`a-dire,
$$
G_{p}(f)\circ \gamma _{p} =\gamma'_{p}\circ G_{p+1}(f),
$$
o\`u $\gamma _{p}:G_{p+1}(X,Y)\longrightarrow G_{p}(X,Y)$ est la
diff\'erentielle du complexe $G_{\bullet}(X,Y)$, qui est une somme
altern\'ee de morphismes de Gysin.
Si $\sigma =(\sigma (1),\cdots,\sigma (p+1))\subseteq I$ et $\nu =(\nu
(1),\cdots,\nu (p))\subseteq J$ sont des suites croissantes, la
composante d'indices $(\sigma ,\nu) $ de
$G_{p}(f)\circ \gamma _{p}$
est, compte tenue de (5.10),
$$ \align
(G_{p(f)}\circ \gamma _{p})_{\sigma, \nu } =&
\sum_{l =1   }^{p+1}
\varepsilon (\sigma  ,\sigma - \sigma (l)  )
m_{\sigma -\sigma ( l)  ,\nu }
f^{*}_{\sigma -\sigma (l)  ,\nu  }i_{\sigma ,\sigma -\sigma (l) *}\\
=&\sum_{l ,\beta   }^{}
\varepsilon (\sigma  ,\sigma - \sigma (l)  )
m_{\sigma - \sigma (l)  ,\nu }
m_{\sigma (l),  \beta }j_{\nu \cup \beta,\nu  *}f^{*}_{\sigma ,\nu
\cup \beta
 },
\endalign
$$
o\`u
$\varepsilon (\sigma ,\mu )=(-1)^{l+1}$
, si  $\sigma  =\mu  \cup \alpha $ est tel que $\sigma  (l)=\alpha
$,  est le signe du morphisme de
Gysin
correspondante \`a l'inclusion $i_{\sigma ,\mu  }:Y_{\sigma
}\longrightarrow Y_{\mu }$. Analoguement, la composante correspondante
de $\gamma'_{p}\circ G_{p+1}(f)$
 est
$$
(\gamma'_{p}\circ G_{p+1}(f))_{\sigma ,\nu }=
\sum_{\beta  }^{}
\varepsilon (\nu \cup \beta  ,\nu)
m_{\sigma  ,\nu \cup \beta   }
j_{\nu  \cup \beta  ,\nu
*}f^{*}_{\sigma  ,\nu \cup \beta }. $$
Alors l'egalit\'e  $$
\sum_{l }^{}
\varepsilon (\sigma  ,\sigma -\sigma (l ) )
m_{\sigma -\sigma (l ) ,\nu }
m_{\sigma (l) , \beta }=
\varepsilon (\nu \cup \beta  ,\nu)
m_{\sigma  ,\nu \cup \beta   }
$$
r\'esulte de la r\`egle de Laplace du d\'eveloppement du d\'eterminant
$m_{\sigma ,\nu \cup\beta }$ par la colonne d'index $\beta $. Cette
fois-ci c'est la fonctorialit\'e du complexe de Koszul qu'on entrevoie
darri\`ere les calculs.
\enddemo

Dor\'enavant nous noterons $G_{\bullet}$ simplement par $G$.

\proclaim{(5.12)\, Proposition} Il existe un foncteur contravariant $$
G:(\Phi ^{op}, \bold{P})\longrightarrow (\Phi
,C^{b}(\Cal{M}^{+}_{rat}(k))) $$ tel que
$G(X,Y)=\bold{s}(h_{*}(S_{\bullet}(Y)\longrightarrow X))(-dim
\, X). $
Le foncteur $G$ v\'erifie les propriet\'es (F1), (F2) et (F3) pour
les \'eclatements avec un centre lisse transverse aux
intersections des  composantes du diviseur $Y$.
 \endproclaim

\demo{Preuve}
La fonctorialit\'e de $G$ est une cons\'equence de la proposition
pr\'ec\'edante et les
propri\'et\'es (F1) et (F2) sont ais\'ement v\'erifie\'es. Prouvons la
propri\'et\'e
(F3). Soient  $(X,Y) $ un
objet de $\bold{P}$, $Z$ un sous-sch\'ema ferm\'e de $X$,
qui est lisse et  transverse \`a chaque intersection $Y_{\sigma }$ des
composantes de $Y$, et $f:X'\longrightarrow
X$ l'\'eclatement de $X$, le long de $Z$. Alors on consid\`ere le
diagramme cart\'esien de $\bold{P}$ d\'efinit par $f$, $(X,Y)$ et
$Z$, $$
 \CD
 (Z',T') @>j>>  (X',Y') \\
@V g VV @VV f V\\
(Z,T) @>i>> (X,Y)
\endCD
$$
Ici $T=Z\cap Y$ est un diviseur de $Z$ de composantes $T_{\alpha
}=Z\cap Y_{\alpha }$, non n\'ecessairement diff\'erentes deux \`a deux.
Ainsi, nous consid\'erons le complexe  $ G (Z,\{T_{\alpha
}\}_{\alpha \in I })$, obtenu avec la famille des $T_{\alpha }$
posiblement repet\'es, et qui est homotopiquenment equivalente \`a
$G(Z,T)$.
 Pour
tout sous-ensemble $\sigma \subseteq
I$, $Z$ est transverse \`a l'intersection $Y_{\sigma }$ , donc la
restriction du diagramme pr\'ec\'edent \`a $Y_{\sigma }$ induit un diagramme
cart\'esien
 $$
 \CD
 T'_{\sigma } @>j_{\sigma }>>  Y'_{\sigma } \\
@V g_{\sigma } VV @VV f_{\sigma } V\\
T_{\sigma } @>i_{\sigma }>> Y_{\sigma }
\endCD
$$
lequel est l'\'eclatement de $Y_{\sigma }$ le long de $T_{\sigma }$.
Donc, on a le diagramme commutatif
$$
 \CD
h(T'_{\sigma }) @<j_{\sigma }^{*}<< h( Y'_{\sigma }) \\
@A g_{\sigma }^{*} AA @AA f_{\sigma }^{*} A\\
h(T_{\sigma }) @<i_{\sigma }^{*}<<h(Y_{\sigma })
\endCD
$$
tous les exc\'es \'etant triviales par transversalit\'e, et ce diagramme  est
acyclique, par (5.3).
Il en r\'esulte que le diagramme
$$
 \CD
G (Z',
\{T'_{\alpha
}\}_{\alpha\in I }
) @<G(j)<< G (X',Y') \\
@A G(g) AA @AA G(f) A\\
 G (Z,
\{T_{\alpha
}\}_{\alpha \in I}
) @<G(i)<<G(X,Y)
\endCD
$$
est acyclique.
\enddemo

 \proclaim{(5.13) \, Th\'eor\`eme}
Il existe un unique foncteur
contravariant
$$
\bold{Sch}(k)\longrightarrow  Ho \,
C^{b}(\Cal{M}^{+}_{rat}(k)), $$
not\'e encore $h$,
tel que:
\roster
\item Si $X$ est un sch\'ema projectif et lisse
sur $k$, $h(X)$ est le motif de Chow associ\'e \`a $X$, $h(X)=(X,id_{X})$.
\item $h$ v\'erifie la propri\'et\'e de descente {\rm (D)}.
\item
Si $X$  est un sch\'ema  projectif et lisse,et
 $Y$ est un diviseur \`a croisements normaux dans $X$ qui est
r\'eunion de diviseurs lisses,
on a un isomorphisme
$h(X-Y)\cong \bold{s}(h_{*}(S_{\bullet}(Y)\longrightarrow X))(-dim\,X)$.
\endroster
\endproclaim
\demo{Preuve}
L'existence de l'extension $h$
et les propri\'et\'es
(F1), (F2), et (F3) sont une cons\'equence de (2.2) et le
lemme suivant.
\enddemo
\proclaim{Lemme 1}
Le foncteur $$
G:(\Phi^{op}, \bold{P})\longrightarrow Ho (\Phi
,C^{b}(\Cal{M}^{+}_{rat}(k))) $$
factorise par le foncteur $\gamma :\bold{P}\longrightarrow
\bold{Sch_{Reg}}(k)$ d\'efini par $\gamma (X,Y)=X-Y$, et il d\'efinit un
foncteur
$$
G:(\Phi ^{op},\bold{Sch_{Reg}}(k))\longrightarrow Ho\,(\Phi ,
C^{b}(\Cal{M}^{+}_{rat}(k))), $$
qui commute avec les foncteurs {\rm tot} et {\rm type}
et  qui v\'erifie (F1), (F2) et (F3).
\endproclaim

\demo{Preuve du lemme 1 }
En effet, d'apr\'es
le th\'eor\`eme de compactification de Nagata et \cite{HC}(I.4), donn\'ee
 une
cat\'egorie ordonnable et finie $\bold{I}$, et $U_{\bullet}$  un
$\bold{I}$-diagramme de $k$-sch\'emas lisses,
il existe une compactification
$U_{\bullet}\longrightarrow X_{\bullet}$
de $U_{\bullet}$ telle que le complementaire
$Y_{i}=X_{i}-U_{i}$ soit un
diviseur \`a croisements normaux dans $X_{i}$, pour tout $i\in \bold{I}$.
Ainsi le foncteur $\gamma $
est essentiellement exhaustif. On montre de la m\^eme fa\c con qu'il est
pleine. Donc, il suffit de montrer que, donn\'e $U_{\bullet}$,
le complexe
$G(X_{\bullet},Y_{\bullet})$ ne d\'epend pas de la compactification
choisie, \`a isomorphisme canonique pr\`es, et m\^eme pour les morphismes.
D'apr\`es
\cite{HC}(I.4) on peut supposer que, donn\'ees deux compatifications
$X_{\bullet 1}$ et $X_{\bullet 2}$ de $U_{\bullet}$ il existe
un morphisme $X_{\bullet 2}\longrightarrow X_{\bullet 1}$. Maintenant on
se ram\`ene au cas o\`u $\bold{I}$ est r\'eduit \`a un point. Mais, dans ce
cas, le morphisme peut \^etre
domin\'e par un \'eclatement de $X_{1}$ avec un centre contenu dans $Y_{1}$.
En raissonnant comme dans (2.3.1), on se ram\`ene au cas o\`u le
morphisme $X_{ 2}\longrightarrow X_{ 1}$
est un \'eclatement de
$X_{1}$ avec un centre lisse contenu dans une intersection de
composantes de $Y_{1}$ et transverse \`a tous les intersections que ne le
contient pas. Alors l'independance r\'esulte du lemme 2  suivante.
Les cas des morphismes s'obtient
de la m\^eme fa\c con en utilisant le diagramme total associ\'e \`a un
morphisme. \enddemo
\proclaim{Lemme 2} Soient $X$ une compactification d'un $k$-sch\'ema
lisse
$U$ telle que $Y=X-U$ est un diviseur  \`a croisements normaux de $X$,
$f:X'\longrightarrow X$ un \'eclatement de $X$ avec un centre lisse $S$
contenu dans une intersection $Y_{\sigma }$ de composantes de $Y$ et
transverse \`a tous les intersections qui ne le contient pas, et
$Y'=f^{-1}(Y)$. Alors, le morphisme $$
G(f):G(X,Y)\longrightarrow G(X',Y')
$$
est un homotopisme.
\endproclaim
\demo{Preuve du lemme 2}
Soient $\widetilde Y$ la transform\'ee stricte de $Y$, $E=f^{-1}(S)$ le
diviseur exceptionnel, $T_{\alpha }=S\cap Y_{\alpha }$, $D_{\alpha
}=E\cap \widetilde Y _{\alpha }$, $\alpha \in I$.
Pour avoir une id\'ee plus claire de la preuve nous consid\'erons d'abord
le cas o\`u $Y$ n'a qu'une seule composante. Dans ce cas on a
$S=T$, et on consid\`ere les complexes de Gysin
$G(S,T)=\bold{s}(id_{*}:h_{*}(T)\longrightarrow h_{*}(S))(dim \,S)$,
qui est acyclique, et $G(E,D)= \bold{s}(i_{D,E*}:h_{*}(D)\longrightarrow
h_{*}(E))(dim \,E)$.

Soit
$\widetilde G (X,Y)$
le complexe simple du diagramme
$$
 \CD
h(S)(-e) @>i_{S,X * }>>h(X) \\
@A i_{*}  AA @AA i_{Y,X*} A\\
h(T)(-e) @>i_{T,Y*}>>h(Y)(-1).
\endCD
$$
o\`u $e$ est la codimension de $S$.
Alors on a une suite exacte scind\'ee
$$0\longrightarrow G(X,Y)@>\varphi >> \widetilde
G(X,Y)\longrightarrow G(S, T)[1]\longrightarrow 0.$$
Comme  $G(S,T)$ est
contractile,  $\varphi $ est un homotopisme.
On a aussi la suite exacte scind\'ee $$
0\longrightarrow G(X,S)\longrightarrow \widetilde G(X,Y)\longrightarrow
G(Y,T)[1]\longrightarrow 0
{}.
$$
D'autre part,
$G(X',Y')$ est
le complexe simple du diagramme
$$
 \CD
h(E)(-1) @>i_{E,X' * }>>h(X') \\
@A i_{*}  AA @AA i_{\widetilde Y,X*} A\\
h(D)(-2) @>i_{D,\widetilde Y*}>>h(\widetilde Y)(-1).
\endCD
$$
Alors on a une suite exacte scind\'ee
$$
0\longrightarrow G(X',E)\longrightarrow
G(X',Y')\longrightarrow G(\widetilde Y,D)[1]\longrightarrow 0
{}.
$$
Le morphisme $G(f):G(X,Y)\longrightarrow G(X',Y') $ factorise par
$$
G(X,Y)@>\varphi >>\widetilde G(X,Y) @>\psi >>G(X',Y'),
$$
o\`u $\psi $ est un morphisme de complexes
qui induit un morphisme des suites exactes pr\'ecedentes,
$$
\CD
0@>>>G(X,S)@>>> \widetilde G(X,Y) @>>> G(Y,T)[1]@>>> 0 \\
@. @V\psi _{0}VV @V\psi VV @V\psi _{1}VV @. \\
0@>>> G(X',E) @>>>G(X',Y')@>>> G(\widetilde Y,D)[1]@>>> 0
\endCD
$$
o\`u $\psi _{0}$ est une homotopisme, d'apr\`es
la suite exacte classique d'un \'eclatement (\cite{M}(\S 9))
appliqu\'ee \`a l'\'eclatement
$$
 \CD
 E @>>> X' \\
@VVV @VVV\\
S @>>> X.
\endCD
$$
et, analoguement, $\psi _{1}$ est un homotopisme, car
$$
 \CD
 D @>>> \widetilde Y \\
@VVV @VVV\\
T @>>> Y
\endCD
$$
est aussi un \'eclatement.
Ceci prouve que  $\psi $ est un homotopisme (cf. preuve de (1.8.3)).

Dans le cas g\'en\'eral, on a des complexes de
Gysin associ\'es $ G(S,\{T_{\alpha} \}_{\alpha })$, et
$G(E,\{D_{\alpha }\}_{\alpha })$, et qui nous noterons simplement
$G(S,T)$ et $G(E,D)$ respectivement.

Soit
$$
\widetilde G (X,Y)=\bold{s}(i_{S,X*}:G(S,T)
\rightarrow G(X,Y)), $$
o\`u $i_{S,X}:(S,T)\longrightarrow (X,Y)$
est le morphisme d'inclusion. Alors l'inclusion canonique
$\varphi :G(X,Y)\longrightarrow \widetilde G(X,Y)$, a une quotient qui
est $G(S, T)[1]
$. Comme $S=\bigcup T_{\alpha }$,  $G(S,T)$ est
contractile, donc $\varphi $ est un homotopisme. On d\'efinit une
filtration croissante finie $\{F_{p}\}_{p\ge 0}$ de  $\widetilde G(X,Y)$
par
 $$
F_{p}(\widetilde G(X,Y)): =
\bold{s}(G_{\le p}(S,T)\longrightarrow G_{\le p}(X,Y)).
$$
Il est \'evidente que $F_{p}$ est un sous-complexe de $\widetilde G
(X,Y)$, et que le gradu\'e par cette filtration  existe et v\'erifie
$$
Gr_{p}^{F}\widetilde G (X,Y)=\bold{s}(G_{p}(S,T)\longrightarrow G_{p}(X,Y)).
$$
D'autre part, $G(X',Y')$ a une d\'ecomposition analogue  \`a celle de
$\widetilde G(X,Y)$ de la
forme $$G(X',Y')=\bold{s}
(i_{E,X'*}:G(E,D)\rightarrow G(X',\widetilde Y)),
$$
qui induit aussi une filtration $F'$ sur $G(X',Y')$ telle que
$$
Gr_{p}^{F'} G (X',Y')=
\bold{s}(G_{p}(E,D)\rightarrow G_{p}(X',\widetilde Y)).
$$
Le morphisme $G(f) $ factorise par
$$
G(X,Y)@>\varphi >>\widetilde G(X,Y) @>\psi >>G(X',Y'),
$$
o\`u $\varphi $ est l'inclusion et $\psi $ est un morphisme de complexes
d\'etermin\'e par les composantes suivantes:
\roster
\item "i)" Le morphisme
$$
G_{p}(f_{S,E})\cup \xi :G_{p}(S,T)\longrightarrow G_{p}(E,D),
$$
o\`u $f_{S,E}:(S,T)\longrightarrow (E,D)$ est la restriction de $f$, et
$\xi $ est la classe de Chern de dimension maximum du fibr\'e d'exc\'es sur
chaque composante, et qui co{\"\i}ncide avec la restriction de la
correspondante classe du diagramme
$$
 \CD
 E @>>> X' \\
@VVV @VVV\\
S @>>> X.
\endCD
$$
\item "ii)" Le morphisme
$$ pr_{1}\circ G_{p}(f):G_{p}(X,Y)\longrightarrow G_{p}(X',\widetilde
Y), $$
o\`u $pr_{1}:G(X',Y')\longrightarrow G(X',\widetilde Y)$ est la projection
canonique.
\item "iii)" Le morphisme
$$
pr_{2}\circ G_{p}(f):G_{p}(X,Y)\longrightarrow G_{p-1}(E,D), p>0,
$$
o\`u $pr_{2}:G_{p}(X',Y')\longrightarrow G_{p-1}(E,D)$ est la projection
canonique.
\endroster
Le morphisme $\psi $ est un morphisme filtr\'e, $\psi (F_{p})\subseteq
F'_{p}$, et la composante $
pr_{2}\circ G_{p}(f):G_{p}(X,Y)\longrightarrow G_{p-1}(E,D)$
n'intervient pas dans le gradu\'e. Ainsi  $$ Gr_{p}(\psi ):
\bold{s}(G_{p}(S,T)\rightarrow G_{p}(X,Y))\longrightarrow
\bold{s}(G_{p}(E, D)\rightarrow G_{p}(X',\widetilde Y))
$$
est d\'efini par les diagrammes
$$
 \CD
h(D_{\sigma  })(-1) @>j_{D,Y * }>>
h(\widetilde Y_{\sigma }) \\ @A f^{*}\cup \xi  AA @AA f^{*} A\\
h(T_{\sigma })(-e) @>i_{T,Y*}>>h(Y_{\sigma })
\endCD
$$
$\sigma \subseteq I$. Puisque ces diagrammes sont contractiles, d'apr\`es
(\cite{M}(\S 9)) \`a nouveau appliqu\'e \`a
l'\'eclatement, possiblement trivial,
 $$
 \CD
 D_{\sigma } @>>> \widetilde Y_{\sigma } \\
@VVV @VVV\\
T_{\sigma } @>>> Y_{\sigma },
\endCD
$$
$\psi $ est un homotopisme filtr\'e, et donc un
homotopisme, par r\'ecurrence sur la longueur de la filtration.
\enddemo
\noindent $\bold{(5.14)}$ Si on applique un foncteur de r\'ealisation,
par exemple le foncteur de cohomologie singuli\`ere, au complexe
$h(X)$ on n'obtient pas la cohomologie singuli\`ere de $X$, mais le terme
$E_{1}$ de la suite spectrale associ\'e \`a la filtration par le poids.
Neanmoins le foncteur $h$ v\'erifie d'autres propri\'et\'es des th\'eories
cohomologiques, par exemple il est ais\'e de voir qu'il v\'erifie la
propri\'et\'e d'homotopie: Si $X$ est
un $k$-sch\'ema et $E\longrightarrow X$ est un fibr\'e vectoriel sur $X$,
$h(X)\longrightarrow h(E)$ est un homotopisme.

 \noindent $\bold{(5.15)}$ Les foncteurs de torsion
$?\otimes
\bold{L}$ (not\'e aussi $?(-1)$)  et de passage au dual, $?^{\wedge}$,
sont des foncteurs exactes dans la cat\'egorie des motifs, donc ils
induissent des morphismes sur le groupe $K_{0}$. Alors, de (5.6),
(5.7),
(5.13) et la dualit\'e de Poincar\'e il r\'esulte imm\'ediatement

\proclaim{Corollaire} Il
existe une application $$
\chi:Ob \,  \bold{Sch}(k) \longrightarrow
{K}_{0}(\Cal{M}^{+}_{rat}(k)) $$
qui v\'erifie:
\roster
\item $\chi(X)=[X]$, si $X$ est un sch\'ema projectif et lisse.
\item $\chi(X)=\chi(\widetilde X)+\chi(Y)-\chi
(\widetilde Y)$, si
$(\widetilde X, \widetilde Y)\longrightarrow (X,Y)$ est un isomorphisme
relatif propre.
\item Si $Y$ est un diviseur lisse d'un sch\'ema lisse
$X$, $$
\chi (X-Y)=\chi  (X)-\chi  (Y)(-1).
$$
\item
Si $X$ est une vari\'et\'e
lisse,
 $\chi
(X)^{\vee }=\chi_{c} (X)(dim \, X)$.
\endroster
En outre, l'application $\chi $ est determin\'e par les conditions {\rm
(1)} \`a {\rm (3)}.
\endproclaim

\remark{$\bold{(5.16)}$ Remarques}1. On
peut v\'erifier  que les ant\'erieures caract\'eristiques
d'Euler $\chi _{c}$ et $\chi $ sont compatibles avec les foncteurs de
r\'ealisation, ainsi par
exemple, elles sont compatibles avec les caract\'eristiques
d'Euler des structures de Hodge mixtes
sur la cohomologie  \`a support compact et la cohomologie, respectivement
(\cite{D1}).

2. Bien que en  r\'ealisation de Betti les caract\'eristiques d'Euler
$\chi _{c}$ et $\chi $ co{\"\i}ncident, en g\'en\'eral $\chi _{c}$ et $
\chi $ ne co{\"\i}ncident pas dans $K_{0}(\Cal{M}^{+}_{rat}(k))$. Ainsi, en
r\'eprennant l'exemple
de Serre du loc. cit., soient $Y$ une vari\'et\'e projective et lisse et $X$
le c\^one affine de base $Y$. Alors, on a
$
\chi _{c}(X)=1+\chi (Y)(-1)-\chi (Y)
$
et
$
\chi (X)=1.
$

\end